\documentclass{article}
\usepackage{latexsym,amsfonts,amsmath,amsthm,makeidx,fontenc}
\usepackage{makeidx}
\makeindex
\newtheorem{definition}{Definition}[section]
\newtheorem{theorem}{Theorem}[section]
\newtheorem{lemma}{Lemma}[section]
\newtheorem{corollary}{Corollary}[section]
\newtheorem{proposition}{Proposition}[section]
\newtheorem{remark}{Remark}[section]
\newcommand{\D}{\mathbb D}
\newcommand{\HS}{\mathbb S}
\newcommand{\RN}{\mathbb R^N}

\newcommand{\iy}{\infty}

\newcommand{\s}{\section}

\newcommand{\R}{\mathbb R}

\newcommand{\lab}{\label}
\newcommand{\bt}{\begin{theorem}}
\newcommand{\et}{\end{theorem}}
\newcommand{\bl}{\begin{lemma}}
\newcommand{\el}{\end{lemma}}
\newcommand{\bd}{\begin{definition}}
\newcommand{\ed}{\end{definition}}
\newcommand{\bc}{\begin{corollary}}
\newcommand{\ec}{\end{corollary}}
\newcommand{\bp}{\begin{proof}}
\newcommand{\ep}{\end{proof}}
\newcommand{\bx}{\begin{example}}
\newcommand{\ex}{\end{example}}
\newcommand{\bi}{\begin{exercise}}
\newcommand{\ei}{\end{exercise}}
\newcommand{\bo}{\begin{proposition}}
\newcommand{\eo}{\end{proposition}}
\newcommand{\br}{\begin{remark}}
\newcommand{\er}{\end{remark}}
\newcommand{\be}{\begin{equation}}
\newcommand{\ee}{\end{equation}}
\newcommand{\ba}{\begin{align}}
\newcommand{\ea}{\end{align}}
\newcommand{\bn}{\begin{enumerate}}
\newcommand{\en}{\end{enumerate}}
\newcommand{\bg}{\begin{align*}}
\newcommand{\bcs}{\begin{cases}}
\newcommand{\ecs}{\end{cases}}

\newcommand{\NN}{{\mathbb N}}

\newcommand{\bean}{\begin{eqnarray*}}
\newcommand{\eean}{\end{eqnarray*}}

\newcommand{\intR}[1]{\int_{\R^N}#1\,  dx}

\numberwithin{equation}{section}

\begin{document}

\title{{\bf On    coupled   Schr\"{o}dinger  systems with  double  critical exponents  and  indefinite  weights}\thanks{Supported by NSFC(11025106,  11371212,  11271386). E-mails: zhongxuexiu1989@163.com, wzou@math.tsinghua.edu.cn}}

\date{}
\author{
{\bf X.  Zhong      {\small \& }  W.   Zou}      \\
\footnotesize   Department of Mathematical Sciences, Tsinghua University,\\
\footnotesize  Beijing 100084, P. R. China}

\maketitle

\vskip0.6in

\begin{center}
\begin{minipage}{120mm}
\begin{center}{\bf Abstract}\end{center}
By using    variational methods, we study  the existence of mountain pass solution    to  the following  doubly critical  Schr\"{o}dinger   system:
$$
  \begin{cases}
  -\Delta u-\mu_1\frac{u}{|x|^2}-|u|^{2^{*}-2}u &=h(x)\alpha|u|^{\alpha-2}|v|^\beta u\quad  \rm{in}\; \R^N,\\
  -\Delta v-\mu_2\frac{v}{|x|^2}-|v|^{2^{*}-2}v &= h(x)\beta |u|^{\alpha}|v|^{\beta-2}v\quad  \rm{in}\; \R^N, \\
  \end{cases}
$$
where $\alpha\geq 2,   \beta\geq 2, \alpha+\beta\leq 2^*$;\;  $ \mu_1, \mu_2\in [0,  \frac{(N-2)^2}{4})$.
The weight function $h(x)$ is allowed to  be sign-changing    so  that the nonlinearities include a large class of indefinite weights.
We show  that  the $PS$ condition is satisfied at  higher energy level when $\alpha+\beta=2^*$ and obtain the existence of mountain pass solution.  Besides,  a nonexistence
result of the ground state is given.

 \vskip0.123in

{\it   Key  words:} Doubly critical system, mountain pass solution, Nehari manifold, indefinite weight.

\vskip0.23in

\end{minipage}
\end{center}

\newpage

\noindent{\bf \Large Content:}

\begin{itemize}

\item [1.]Introduction
\item [2.]  Nehari Manifold.
\item [3.]  Analysis of the Palais-Smale Sequences
\item [4.] Nonexistence of the Nontrivial Least Energy Solution
\item [5.] The Existence of Mountain Pass Solutions
\begin{itemize}
\item [5.1]  The case of $N=3, \frac{1}{2}<\frac{1-4\mu_1}{1-4\mu_2}, S^{\frac{N}{2}}(\mu_2)+S^{\frac{N}{2}}(\mu_1)\leq S^{\frac{N}{2}}$
\item [5.2]  The case of $N=4, \frac{1}{2}<\big(\frac{1-\mu_1}{1-\mu_2}\big)^{\frac{3}{2}}, S^{\frac{N}{2}}(\mu_2)+S^{\frac{N}{2}}(\mu_1)\leq S^{\frac{N}{2}}$
\item [5.3] The case of $N=3, \frac{1}{2}<\frac{1-4\mu_1}{1-4\mu_2}, 2\big(\frac{S(\mu_2)+S(\mu_1)}{2}\big)^{\frac{N}{2}}> S^{\frac{N}{2}}$
\item [5.4] The case of $N=4, \frac{1}{2}<\big(\frac{1-\mu_1}{1-\mu_2}\big)^{\frac{3}{2}}, 2\big(\frac{S(\mu_2)+S(\mu_1)}{2}\big)^{\frac{N}{2}}> S^{\frac{N}{2}}$

\end{itemize}
\end{itemize}

\vskip0.6in

\s{Introduction}

\renewcommand{\theequation}{1.\arabic{equation}}
In this paper, we investigate the existence of  solutions  to  the following nonlinear Schr\"{o}dinger  system:
\be\lab{hhe6-4}
\begin{cases}
-\Delta u-\mu_1 \frac{u}{|x|^2}-|u|^{2^*-2}u& = h(x)\alpha |u|^{\alpha-2}|v|^\beta u, \;\;  \hbox{in}\; \R^N,\\
-\Delta v-\mu_2\frac{v}{|x|^2}-|v|^{2^*-2}v & = h(x)\beta |u|^\alpha |v|^{\beta-2}v,  \;\;   \hbox{in}\; \R^N,
\end{cases}
\ee
where $\alpha\geq 2, \beta\geq 2, \alpha+\beta\leq 2^*$; $\mu_1, \mu_2 \in [0,  \Lambda_N), \Lambda_N:= \frac{(N-2)^2}{4}; h(x)\in L^\infty(\R^N)$.
The interest   for such  systems  is motivated by its applications to plasma physics, nonlinear optics, condensed matter physics, etc. For example, the coupled nonlinear Schr\"{o}dinger systems arise in the description of several physical phenomena such as the propagation of pulses in birefringent optical fibers and Kerr-like photorefractive media, see \cite{AN,KIP,CRM1,CRM2,PCH,SJJN},  etc.  Also, it is related to
the following Gross-Pitaevskii equations (cf. \cite{HMEWC, TV}):
\be\label{eq0001}
\begin{cases}
-i\frac{\partial}{\partial t}\Phi_1=\Delta \Phi_1- a(x)\Phi_1+
\mu_1 |\Phi_1|^2 \Phi_1+\nu |\Phi_2|^2\Phi_1,\,\, x\in \RN, \,\,t>0,\\
-i\frac{\partial}{\partial t}\Phi_2=\Delta \Phi_2- b(x)\Phi_2+
\mu_2|\Phi_2|^2 \Phi_2+\nu |\Phi_1|^2\Phi_2,\,\, x\in \RN, \,\,t>0,\\
\Phi_j=\Phi_j(x,t)\in\mathbb{C},\quad j=1,2,\\
\Phi_j(x,t)\to 0,  \quad\hbox{as}\,\, |x|\to+\iy, \,\,t>0, \,\,j=1,2,
\end{cases}\ee
where $i$ is the imaginary unit; $a(x), b(x)$ are potential functions. Problem (\ref{eq0001}) also arises
in  the Hartree-Fock theory for a double condensate, i.e., a binary mixture of Bose-Einstein
condensates in two different hyperfine states (see \cite{EGBB}).

\vskip0.1in

We call a solution $(u, v)$ nontrivial if both $u\not\equiv 0$ and $v\not\equiv 0$;
we call a solution $(u, v)$ semi-trivial if $(u, v)$ is a type of $(u, 0)$ or $(0, v)$. The existence of semi-trivial solution is equivalent to the solution of the following  scalar  equation:
\be\lab{LC-1}
-\Delta u-\mu\frac{u}{|x|^2}=|u|^{2^*-2}u,\;\hbox{in}\;\R^N,
\ee
whose solutions  have  been  figured out.   Here, when  $\mu=0$,  we refer the readers to  \cite{Aubin,AJ,AL,HP,HB,HL,BWL,BWL2}.
 When $\mu\in(0, \frac{(N-2)^2}{4})$, by \cite{T.S}, an explicit solution of (\ref{LC-1}) exists, namely
\be\lab{zzz-2}
z_1^\mu(x):=\frac{A(N,\mu)}{|x|^{a_\mu}\big(1+|x|^{2-\frac{4a_\mu}{N-2}}\big)^{\frac{N-2}{2}}};
\ee
where $  a_\mu:=\frac{N-2}{2}-\sqrt{(\frac{N-2}{2})^2-\mu}$  and $  A(N,\mu):=\Big[\frac{N(N-2-2a_\mu)^2}{N-2}\Big]$.
This solution is also known to be the unique positive solution up to a conformal transformation of the form
\be\lab{zzz-1}
z_\sigma^\mu=\sigma^{-\frac{N-2}{2}}z_1^\mu(\frac{x}{\sigma}), \;\;\; \sigma>0.
\ee
Before returning to  the existence  and nonexistence of the nontrivial solutions of  (\ref{hhe6-4}),   we recall  the very recent paper  \cite{BVI}, where the authors studied  the existence of   solutions to the following system:
\be\lab{hhe6-1}
\begin{cases}
-\Delta u-\mu_1 \frac{u}{|x|^2}-|u|^{2^*-2}u & =\nu \cdot  h(x)\alpha |u|^{\alpha-2}|v|^\beta u,\;\hbox{in}\; \R^N,\\
-\Delta v-\mu_2\frac{v}{|x|^2}-|v|^{2^*-2}v & =\nu  \cdot h(x)\beta |u|^\alpha |v|^{\beta-2}v,\;\hbox{in}\; \R^N,
\end{cases}
\ee
where $\mu_1, \mu_2\in (0, \frac{(N-2)^2}{4})$ and the parameter $\nu$ serves  as a regulator. Throughout the paper \cite{BVI},   $h(x)$ satisfies  the following condition:
\be\lab{zwm1} h(x)\geq 0, h(x)\not\equiv 0, h(x)\in L^1(\R^N)\cap L^\infty(\R^N).\ee  Let
\be\label{zwm2}S(\mu_i):= \left(1-\frac{4\mu_i}{(N-2)^2}\right)^{\frac{N-1}{N}}S,\ee
where  $S$ is the sharp constant of $D^{1,2}(\RN)\hookrightarrow L^{2^\ast}(\RN)$  satisfying
\be\label{sobolev}\intR{|\nabla u|^2}\ge
S\left(\intR{|u|^{2^\ast}}\right)^{\frac{2}{2^\ast}}\ee
and
\be\label{sobolev-2}
S=\frac{N(N-2)}{4}|\HS^N|^{\frac{2}{N}}=\frac{N(N-2)}{4}2^{\frac{2}{N}}\pi^{1+\frac{1}{N}}\Gamma\big(\frac{N+1}{2}\big)^{-\frac{2}{N}}.
\ee
In particular,  $S=3(\frac{\pi}{2})^{\frac{4}{3}}$ when $\displaystyle N=3$; $ S=\frac{4\sqrt{6}\pi}{3}$ when $\displaystyle N=4$.
\vskip0.3in

If  the coupling terms are of subcritical case, i.e.,  $\alpha+\beta<2^*$,  when $\max\{\alpha, \beta\}<2$ the authors of \cite{BVI} prove  that the least  energy $c$ satisfies  $$c< \frac{1}{N}\big(\min\{S(\mu_1), S(\mu_2)\}\big)^{\frac{N}{2}}$$ and obtain the existence of nontrivial ground state solution;  when  $\max\{\alpha, \beta\}=2$, the similar results hold provided  that the regulator $\nu$ is large enough. If $\min\{\alpha, \beta\}>2$, the ground state energy is achieved by and only by semi-trivial solutions; if $\min\{\alpha, \beta\}=2$, the similar results hold provided the regulator $\mu$ small enough. When $N=3$ and $S^{\frac{N}{2}}(\mu_1)+S^{\frac{N}{2}}(\mu_2)<S^{\frac{N}{2}}, \frac{\Lambda_N-\mu_1}{\Lambda_N-\mu_2}>\frac{1}{2}$, they obtain  the existence of mountain pass solution provided that $\nu$ is sufficiently small.

\vskip0.23in

For the critical case, that is,  $\alpha+\beta=2^*$,  \cite{BVI} assumed  that   $h(x)$ is a radial function satisfying
$$\begin{cases}
h\in L^\infty(\R^N), h\geq 0, h\not\equiv 0, h\;\hbox{is continuous in a neighborhood of $0$ and $\infty$},\\
h(0)=\underset{|x|\rightarrow \infty}{\lim}h(x)=0.
\end{cases}$$
Then they   obtained  that if $N\geq 5, \max\{\alpha, \beta\}<2$, then (\ref{hhe6-1}) possesses a nontrivial ground state solution; if $N=3,4, \min\{\alpha, \beta\}\geq 2, \mu_2\leq \mu_1<\frac{(N-2)^2}{4}, \frac{\Lambda_N-\mu_1}{\Lambda_N-\mu_2}>2^{-\frac{2}{N-1}}$, they obtain the existence of mountain pass solution provided that the regulator $\nu$ is small enough. But for the case of $\alpha+\beta=2^*$ and  $h(x)$ is not radial, they only obtain the existence of ground state solution for $\max\{\alpha, \beta\}<2$ and $\nu$ small enough .

\vskip0.1in

Note that if  $1<\alpha<2, 1<\beta<2, \alpha+\beta=2^*$,  hence $2^*<4$,  which means that the results in \cite{BVI} do not include the case of dimension $N=3, 4$ if $h(x)$ is not radial. If $\mu_2<\mu_1$ and $\beta<2$, they also obtain the existence of ground state solution provided $\nu$ small enough.  However,  if $\mu_1=\mu_2$ or $\min\{\alpha, \beta\}\geq 2$, whether there still exists a nontrivial solution remains  open.  At the end of \cite{BVI}, the authors also impose a list of complicated conditions on $h(x)$ and emphasize that if $h(x)$ has a fixed sign, by using  the perturbation argument, they obtain the existence of nontrivial solutions provided that $\nu$ is small enough. In \cite[Theorem 3.8]{BVI}, they considered the case of $\alpha+\beta<2^*, \alpha\geq 2, \beta\geq 2$. Note that for this case   there must hold $N=3$.  More important, in order to prove the Palais-Smale compactness condition, they need that
\be\lab{zwm21}S^{\frac{N}{2}}(\mu_1)+S^{\frac{N}{2}}(\mu_2)<S^{\frac{N}{2}}.\ee
Under this hypothesis,  they obtained a mountain pass solution provided that $\nu$ is small enough.   We emphasize  that (\ref{zwm21})
can not hold for many  ranges of $\mu_1$ and $ \mu_2$; for example:  $\mu_1=0$ or $ \mu_2=0$, or $\mu_1,\mu_2>0 \;\hbox{such that}\;\mu_1+\mu_2\leq\frac{1}{4}$.

\vskip0.16in

Naturally, we concern the following  questions which are still standing   open before us.
\begin{itemize}
\item[1)] Whether the role of the parameter $\nu$ is essential and can be dropped?
\item[2)] What happens if  $h(x)$ is not radial?
\item[3)] What happens if  $h(x)$ is sign-changing?
\item[4)] What happens if   (\ref{zwm21}) is not true, i.e., $S^{\frac{N}{2}}(\mu_1)+S^{\frac{N}{2}}(\mu_2)\geq S^{\frac{N}{2}}$?
\item[5)] When $\alpha+\beta=2^*$, whether there exists  a  mountain pass solution  to (\ref{hhe6-4})?
\end{itemize}
The main purpose of the present paper is to study the existence of mountain pass solution when the ground state energy is only achieved by semi-trivial solution. We will always assume $h(x)$ is sign-changing and not necessary radial.  Now  it is the place to state our results in the current paper.   We need one of the following two conditions:
\begin{itemize}
\item[(${\bf H_1}$)] $h(x)\in L^{\frac{N}{2}}(\R^N)$;
\item[(${\bf H'_1}$)]$h(x)$ is continuous in $\R^N$,\;\;  $h(0)\leq 0,\; \underset{|x|\rightarrow \infty}{\limsup}h(x)\leq 0$.  Moreover,  we assume    that $\gamma:=  \|h_{-}\|_\infty \max\{\alpha, \beta\}<1$, where $h_{-}:=\min\{h, 0\}$.
\end{itemize}
Further, we suppose that  $h(x)$ satisfies the following integrable  condition:
\begin{itemize}
\item[(${\bf H_2}$)] $\int_{\R^N}h(x)|z_{\sigma}^{\mu_1}|^\alpha |z_{\sigma}^{\mu_2}|^\beta dx > 0$ for some $\sigma>0$, where $z_\sigma^\mu$ is defined in (\ref{zzz-1}).
\end{itemize}

 Let us denote
\be\lab{ha1}\Theta :=\begin{cases}\|h_+(x)\|_{L^{\frac{2^*}{2^*-\alpha-\beta}}(\R^N)}\quad &\hbox{if}\; \alpha+\beta<2^*,\\
\|h_+(x)\|_{L^{\infty}(\R^N)}\quad & \hbox{if}\; \alpha+\beta=2^*,  \end{cases}\ee
where $h_+:=\max\{h, 0\}$. Without loss of generality, throughout this paper we always assume  $\mu_2\leq \mu_1$.

\subsection{Nonexistence of the Nontrivial Least Energy Solution}

 The first main result of the current paper   concerns with the nonexistence of the ground state to   (\ref{hhe6-4}) for all $N\geq 3$ .

\bt\lab{h-th5-1-1}
Assume that  either $  \beta\geq 2, \mu_2<\mu_1$ or $\alpha\geq 2, \beta\geq 2, \mu_2=\mu_1=\mu$.  Further, suppose that  $$ \begin{cases}&h(x)\;\hbox{satisfies}\; (H_1)\;\hbox{  if  }\;\alpha+\beta<2^*;\\   &h(x)\;\hbox{satisfies}\;(H'_1)   \hbox{  if  }\;\alpha+\beta=2^*. \end{cases}$$
Then there exists $\Theta_0>0,$  depending on $N, \alpha, \beta, \mu_1, \mu_2$, such that  if $\Theta\leq \Theta_0$, then the least energy
of the system is exactly  equal to $\frac{1}{N}S^{\frac{N}{2}}(\mu_1).$
Moreover, it is achieved by and only by

\begin{itemize}
\item   $(\pm z_{\sigma}^{\mu_1}, 0)$    if $\mu_2<\mu_1$, where  $\sigma>0$;
\item   $(\pm z_\sigma^\mu, 0)$ and $(0, \pm z_\sigma^\mu)$ if $\mu_2=\mu_1=\mu\neq0$;
\item   $(\pm z_{\sigma,x_i}, 0)$ and $(0, \pm z_{\sigma,x_i})$ if $\mu_2=\mu_1=0$,  where
$$z_{\sigma, x_i}=\sigma^{-\frac{N-2}{2}}z_1(\frac{x-x_i}{\sigma}),  z_1(x)=\frac{[N(N-2)]^{\frac{N-2}{4}}}{[1+|x|^2]^{\frac{N-2}{2}}}, \quad \sigma>0, x_i\in\R^N. $$
 \end{itemize}
That is,  problem (\ref{hhe6-4}) has no nontrivial least energy solution.
\et

\br\lab{h-remark5-1}  In the above theorem,  the constant $\Theta_0$ has an explicit formula in terms of  $N, \alpha, \beta, \mu_1, \mu_2$. To avoid tedious notations, we prefer to
give them in Section 4.  For the system (\ref{hhe6-1}), the authors of \cite{BVI} had constructed similar results (see\cite[Theorem 3.4]{BVI}). But they required that $\alpha, \beta\geq 2$ and   $h(x)\geq 0$.
Here we  improve  the results of \cite[Theorem 3.4]{BVI} to the system (\ref{hhe6-4}). When $\mu_1\neq \mu_2$, we only  require $ \beta\geq2$. Moreover,   $h(x)$  is allowed to be   sign-changing in our case.
\er

\newpage

\subsection{Mountain pass solution: the case of   $N=3.$}

 In this case, $ \Lambda_N=\frac{1}{4}.$

\bt\lab{th-ling1}Assume $N=3, \alpha\geq2,\beta\geq2, \alpha+\beta<2^*$ and   $h(x)$ satisfies $(H_1)$ and $(H_2)$.  Furthermore,
assume that $ \frac{1}{2}<\frac{1-4\mu_1}{1-4\mu_2}$ and that
\be\lab{zwm918918}\displaystyle \hbox{ either }\;\;   S^{\frac{N}{2}}(\mu_2)+S^{\frac{N}{2}}(\mu_1)\leq S^{\frac{N}{2}}\;\;\hbox{  or } \;\; \displaystyle 2\Big(\frac{S(\mu_1)+S(\mu_2)}{2}\Big)^{\frac{N}{2}}>S^{\frac{N}{2}}.\ee   Assume further
 \be\lab{d111}\Theta\leq  10^{-4}\big[(2-8\mu_1)^{\frac{2}{3}}-\big(1-4\mu_2\big)^{\frac{2}{3}}\big],\ee
then the  problem (\ref{hhe6-4}) has a nontrivial weak solution $(u_0, v_0)$ with  $u_0\geq 0, v_0\geq0, u_0v_0\not\equiv 0$.
\et

 \br\lab{zwm918-2012}If $\mu_2=\mu_1$, then the alternatives of (\ref{zwm918918}) hold true automatically.\er

\vskip0.1336in
The next theorem is about the   case of  $N=3$ and $ \alpha+\beta=2^*$, which means that the coupling terms are of  critical.
\bt\lab{th-ling2}
Assume $N=3, \alpha\geq2,\beta\geq2, \alpha+\beta=2^*$ and     $h(x)$ satisfies $(H'_1)$ and $(H_2)$. Furthermore,
assume that $\mu_1+\mu_2\neq 0, \frac{1}{2}<\frac{1-4\mu_1}{1-4\mu_2}$ such that  either \be\lab{wmz918} S^{\frac{N}{2}}(\mu_2)+S^{\frac{N}{2}}(\mu_1)\leq S^{\frac{N}{2}}\;\;\;\hbox{or}\;\;\;   2\Big(\frac{S(\mu_1)+S(\mu_2)}{2}\Big)^{\frac{N}{2}}>S^{\frac{N}{2}}.\ee
Assume  further
\be\lab{d222}
\Theta\leq \min\Big\{\frac{\mu_1+\mu_2}{12}, 10^{-4}\big[(2-8\mu_1)^{\frac{2}{3}}-\big(1-4\mu_2\big)^{\frac{2}{3}}\big]\Big\},\ee
then   the problem (\ref{hhe6-4}) has a nontrivial weak solution $(u_0, v_0)$ such that $u_0\geq 0, v_0\geq0, u_0v_0\not\equiv 0$.
\et

\vskip0.336in

\subsection{Mountain pass solution:  the case of   $N=4.$}

  For the case of $N=4$, we know  $2^\ast=4$  and $ \Lambda_N=1.$ If  $\alpha\geq2,  \beta\geq2, \alpha+\beta=2^*$, we must have
$ \alpha=\beta=2$. Thus,  (\ref{hhe6-4}) becomes  a type of Bose-Einstein Condensates (BEC) equation   in $\R^4$:
\be\lab{hsn-e1}
\begin{cases}
-\Delta u-\mu_1\frac{u}{|x|^2}&=u^3+2 h(x)v^2u,\;\;\;\hbox{in}\;\R^4,\\
-\Delta v-\mu_2\frac{v}{|x|^2}&=v^3+2 h(x)u^2v,\;\;\;\hbox{in}\;\R^4,\\
u\geq 0, v\geq 0.
\end{cases}
\ee
Note that both the  cubic   terms ($u^3$ and $ v^3$) and the coupling terms ($v^2u$ and $ u^2v$) on the right-hand
sides of   (\ref{hsn-e1}) are of critical growth.

\bt\lab{th-ling3}  Assume   $(H'_1)$ and $(H_2)$.
Suppose  $\mu_1+\mu_2\neq 0, \frac{1}{2}<\big(\frac{1-\mu_1}{1-\mu_2}\big)^{\frac{3}{2}}$ such that
\be\lab{wmz91800}\displaystyle  \hbox{either }\;\; S^{\frac{N}{2}}(\mu_2)+S^{\frac{N}{2}}(\mu_1)\leq S^{\frac{N}{2}}\;\;\; \hbox{ or  }\;\;\;
\displaystyle 2\Big(\frac{S(\mu_1)+S(\mu_2)}{2}\Big)^{\frac{N}{2}}>S^{\frac{N}{2}}.\ee    Assume further
\be\lab{d333}\Theta\leq\frac{2-(1-\mu_1)^{\frac{3}{2}}-(1-\mu_2)^{\frac{3}{2}}}{16},\ee
 then  the problem (\ref{hsn-e1}) has a nontrivial weak solution
 $(u_0, v_0)$ such that $u_0\geq 0, v_0\geq0, u_0v_0\not\equiv 0$.
\et

\br\lab{13zoujia1}  Basically,      the upper bounds  of $\Theta$ in the above theorems are not sharp.  However, in order to
determine an unambiguous  range of $\Theta$,  we prefer to give the explicit formulas for those constants.
The optimal range  of $\Theta$ is an interesting open question. \er

\br\lab{main-difficulties}
One of the main difficulties of studying this kind of problems is the failure of the (PS)
condition due to the critical term $|u|^{2^*-2}u$ and the unbounded domain $\R^N$, especially for the couple terms $|u|^{\alpha-2}u|v|^\beta$ and $|u|^\alpha |v|^{\beta-2}v$ with $\alpha+\beta=2^*$. One has  to overcome the difficulties on determining the compactness threshold. People usually study the case of $c<\frac{1}{N}S^{\frac{N}{2}}$. One of the main innovation of our present work is that we obtain a nontrivial solution with energy higher than $\frac{1}{N}S^{\frac{N}{2}}$.
We also have to overcome the difficulties brought by the indefinite sign of the weight function $h(x)$, especially when $h(x)$ is not radial.
\er


 \vskip0.33in
\s{Nehari Manifold}

\renewcommand{\theequation}{2.\arabic{equation}}
Let ${\D}:=D^{1,2}(\R^N)\times D^{1,2}(\R^N)$,  where  $D^{1,2}(\R^N)$ is the completion of $C_0^\infty(\R^N)$ with respect to the norm
$$\|u\|_{D^{1,2}(\R^N)}:=\Big(\int_{\R^N} |\nabla u(x)|^2 dx\Big)^{1/2}.$$
By the Hardy inequality, when $0<\mu<\frac{(N-2)^2}{4}$, $\|u\|_{D^{1,2}(\R^N)}$ is equivalent to the following norm:
$$\|u\|_{\mu}:=\Big(\int_{\R^N} (|\nabla u(x)|^2-\mu\frac{u^2}{|x|^2}) dx\Big)^{1/2}.$$
For simplicity,  we will also use the notation of $\|u\|_0$ to represent $\|u\|_{D^{1,2}(\R^N)}$.
For $(u, v)\in {\D}$, define  the norm
$$\|(u,v)\|_{\D}=\Big(\|u\|_{\mu_1}^{2}+\|v\|_{\mu_2}^{2}\Big)^{1/2}.$$
A pair of function $(u, v)$ is said to be a weak solution of problem (\ref{hhe6-4}) iff
\begin{align}\lab{zzzz-2}
&\int_{\R^N}\nabla u\cdot \nabla \varphi_1 dx-\mu_1\int_{\R^N}\frac{u\varphi_1}{|x|^2}dx+\int_{\R^N}\nabla v\cdot \nabla \varphi_2 dx-\mu_2\int_{\R^N}\frac{u\varphi_2}{|x|^2}dx\nonumber\\
&-\int_{\R^N}|u|^{2^*}u\varphi_1dx-\int_{\R^N}|v|^{2^*-2}v\varphi_2dx-\alpha\int_{\R^N}h(x)|u|^{\alpha-2}u|v|^\beta\varphi_1dx\nonumber\\
&-\beta\int_{\R^N}h(x)|u|^\alpha|v|^{\beta-2}v\varphi_2dx=0\;\;\;\hbox{for all $(\varphi_1, \varphi_2) \in {\D}$.}
\end{align}
Thus, the corresponding energy functional of problem (\ref{hhe6-4}) is defined by
\begin{equation}\lab{ee1}
\Phi(u, v)=\frac{1}{2}||(u, v)||_{{\D}}^{2}-\frac{1}{2^{*}}\Big(||u||_{L^{2^{*}}(\R^N)}^{2^{*}}+||v||_{L^{2^{*}}(\R^N)}^{2^{*}}\Big)
-\int_{\R^N}h(x)|u|^{\alpha}|v|^{\beta}dx
\end{equation}
for all $(u, v)\in {\D}$. The associated Nehari manifold is defined as
$${\mathcal{N}}:=\Big\{(u, v)\in {\D}\backslash\big\{(0,0)\big\}:J(u, v)=0\Big\},$$
where
\begin{eqnarray}\lab{J}
J(u, v)&=&\Big\langle\Phi'(u, v), (u, v)\Big\rangle=
||(u, v)||_{\D}^{2}-\Big(||u||_{L^{2^{*}}(\R^N)}^{2^{*}}\nonumber\\
&&+||v||_{L^{2^{*}}(\R^N)}^{2^{*}}\Big)-(\alpha+\beta)\int_{\R^N}h(x)|u|^{\alpha}|v|^{\beta}dx,
\end{eqnarray}
 and $\displaystyle \Phi'(u, v)$ denotes the Fr$\acute{e}$chet derivative of $\Phi$ at $(u, v)$, $\displaystyle \langle\cdot, \cdot\rangle$ is the duality product between ${\D}$ and its dual space ${\D}^\ast$.  We  have the following properties on  the Nehari manifold.

\bl\lab{group-lm1} Assume $\alpha+\beta\leq 2^{*}$. In particular,   if $\alpha+\beta=2^*$, we require $(H'_1)$ instead of $(H_1)$.  Then $\forall (u, v)\in {\D}\backslash\big\{(0, 0)\big\}$, there exists a unique $t=t_{(u, v)}>0$ such that $t(u, v)=(tu, tv)\in {\mathcal{N}}$. Furthermore,  there exists $\delta>0$ such that $t_{(u, v)}\geq\delta$ for all $(u, v)\in \mathcal{S}:=\Big\{(u, v)\in {\D}:   ||(u, v)||_{{\D}}^{2}=1\Big\},$ and ${\mathcal{N}}$ is closed and bounded away from $(0, 0)$.
\el
\bp
For $(u, v)\in \D \backslash\big\{(0, 0)\big\}$, we denote that
\begin{equation}\lab{ee4}
\begin{cases}
a:=||u||_{L^{2^{*}}(\R^N)}^{2^{*}}+||v||_{L^{2^{*}}(\R^N)}^{2^{*}}>0;\\
b:=(\alpha+\beta)\int_{\R^N}h(x)|u|^\alpha|v|^\beta dx;\\
c:=||(u, v)||_{\D}^{2}>0.
\end{cases}
\end{equation}Then,
$ \frac{d}{dt}\Phi(tu, tv)=-tg(t),$
where $\displaystyle g(t):=a\;t^{2^{*}-2}+b\;t^{\alpha+\beta-2}- c.$
Firstly we consider the case $b\geq 0$.   Note that
there  exists a unique  $t_0>0$ such that $g(t_0)=0$. Moreover, $g(t)<0$ for $0<t<t_0$ and
$ g(t)>0$  for $t>t_0$.
 Secondly we consider the case $b<0$. If $\alpha+\beta<2^*$,
  there exists some $s>0$ such that $g(s)=0$. Let $t_0$ be the minimum of the solutions of $g(t)=0$, that is,
$g(t_0)=0$  and $g(t)<0$ for $t<t_0$.
For   $\forall \;t>t_0$, it is easy to check that  $g'(t)>0$.
Thus, $g(t)>0$ for $t>t_0$. Then, $t_0$ is the unique solution to $g(t)=0$.
If $\alpha+\beta=2^*$, by the Young's inequality and $(H'_1)$, we have
\begin{eqnarray}\lab{5-1}
|b|&\leq&\big|2^{*}\int_{\R^N}h_{-}(x)|u|^\alpha|v|^\beta dx\big|\nonumber\\
&\leq& 2^{*} \|h_{-}\|_\infty  \big(\frac{\alpha}{2^{*}}||u||_{L^{2^{*}}(\R^N)}^{2^{*}}+\frac{\beta}{2^{*}}||v||_{L^{2^{*}}(\R^N)}^{2^{*}}\big)\nonumber\\
&\leq& \|h_{-}\|_\infty  \max\{\alpha, \beta\}a\nonumber\\
 &<& a.
\end{eqnarray}
Hence, $a+b>0$. Then there exists a unique positive solution $t_0$ to the equation: $g(t)=(a+b)t^{2^{*}-2}-c=0$.
In particular,
$g(t)<0$ for $0<t<t_0$ and
$ g(t)>0$  for $t>t_0$.  Let $t_{(u,v)}:=t_0$ be  defined as above,   we finally obtain that
\begin{displaymath}
\frac{d}{dt}\Phi(tu, tv)=-tg(t)
\begin{cases}
>0 \quad &\text{for $0<t<t_{(u, v)}$};\\
<0 \quad & \text{for $t>t_{(u, v)}$.}
\end{cases}
\end{displaymath}
In either case, there exists a unique $t_{(u, v)}>0$ such that
$\displaystyle \Phi\Big(t_{(u, v)}u, t_{(u, v)}v\Big)=\underset{t>0}{\max}\Phi(tu, tv) $
and  $t_{(u, v)}(u, v)\in {\mathcal{N}}.$
For $\omega=(u, v)\in {\mathcal{S}}$,  since $h(x)\in L^{\frac{N}{2}}(\R^N)\cap L^{\infty}(\R^N)$ and $\alpha+\beta\leq 2^{*}$, there exists some $C>0$ such that
\begin{displaymath}
\begin{cases}
a&=||u||_{L^{2^{*}}(\R^N)}^{2^{*}}+||v||_{L^{2^{*}}(\R^N)}^{2^{*}}\leq C;\\
|b|&=\big|(\alpha+\beta)\int_{\R^N}h(x)|u|^{\alpha}|v|^{\beta}dx\big|\leq C;\\
c&=1.
\end{cases}
\end{displaymath}
We consider the equation $\displaystyle g(t)=a\;t^{2^{*}-2}+b\;t^{\alpha+\beta-2}-1=0$.
If $\displaystyle b\leq 0$, we have $\displaystyle at^{2^{*}-2}\geq 1$, hence $\displaystyle t\geq ({a})^{-\frac{1}{2^{*}-2}}\geq ({C})^{-\frac{1}{2^{*}-2}}$.
If $\displaystyle b>0$, then we have either
$ a t^{2^{*}-2}\geq \frac{1}{2}\quad \text{or}\quad b t^{\alpha+\beta-2}\geq \frac{1}{2}.$
Thus, either
$ t\geq \Big(\frac{1}{2C}\Big)^{\frac{1}{2^{*}-2}}$ or $ t\geq \Big(\frac{1}{2C}\Big)^{\frac{1}{\alpha+\beta-2}}.$
Therefore, there exists $\delta>0$ such that
$ t_{(u, v)}\geq \delta \quad \text{for all }  (u, v)\in S.$
Therefore, ${\mathcal{N}}$ is bounded away from $(0, 0)$. Obviously,   ${\mathcal{N}}$ is closed.
\ep


\vskip0.33in
\s{Analysis of the Palais-Smale Sequences}

\renewcommand{\theequation}{3.\arabic{equation}}
In this  section, we perform a careful analysis of the behavior of the Palais-Smale sequences with the aid of the  concentration-compactness principle in \cite{PL1,PL2}, which allows to recover  the compactness below some critical threshold.   Set

\be\label{zwm-may-3-} \widetilde{C}_{N,\alpha,\beta}:=\frac{\Big(1-\frac{4\max\{\mu_1,\mu_2\}}{(N-2)^2}\Big)^{\frac{1-N}{N}}-1}{\max\{\alpha,\beta\}}.\ee

\vskip0.136in

\bl\lab{lm-6-3}
Assume $\mu_2=0, \mu_1=\mu\in(0, \frac{(N-2)^2}{4})$ and

$$ \hbox{ either } \begin{cases}& (H_1) \\&\alpha+\beta<2^*\end{cases}\;\;  or \;\; \begin{cases}&(H'_1) \\&\alpha+\beta=2^*\end{cases}.$$ Let $\{(u_n,v_n)\}\subset\mathcal{N}$ be a Palais-Smale sequence for $\Phi|_{\mathcal{N}}$ at level $c<\frac{1}{N}S^{\frac{N}{2}}(\mu)$. Then, there exists some constant $C$, such that  $||(u_n,v_n)||_{\D}\leq C$
for all $n\in \NN$ and $\Phi'(u_n, v_n)\rightarrow 0$ in the dual space ${\D}^{*}$. Moreover,
\begin{itemize}
\item[(1)]for the case of $\begin{cases}&(H_1) \\ &\alpha+\beta<2^*\end{cases}$, we have  $(u_n,v_n)\rightarrow (u_0, v_0)$ in $\D$  up to a subsequence;\\
\item[(2)]for the case of $\begin{cases}&(H'_1) \\&\alpha+\beta=2^*\end{cases}$, if $h(x)$ is radial, we have   $(u_n, v_n)\rightarrow (u_0, v_0)$ in ${\D}$   up to a subsequence.  However, if $h(x)$ is not radial, we   obtain the same result provided  that the additional hypothesis $ \Theta<  \widetilde{C}_{N,\alpha,\beta} $  (see   (\ref{zwm-may-3-}))  holds.
\end{itemize}

\el
\bp  The ideas for proving this lemma are quite similar to the cases of $\mu_1,\mu_2>0$  and $\mu_1=\mu_2=0$ in \cite[Lemma 6.3]{ZZ},  we  omit the details.
\ep
\vskip 0.236in

 In the next section, we will study the nonexistence of nontrivial ground state solutions.  In view of  the nonexistence of  the ground state to the system,    we will  investigate the existence of  mountain pass solutions of system (\ref{hhe6-4}).  For this goal, we need an improved Palais-Smail condition at  higher energy level.   Let  \be\lab{hhe6-2}
I_\mu(u)=\frac{1}{2}\int_{\R^N}\big(|\nabla u|^2-\mu\frac{u^2}{|x|^2}\big) dx-\frac{1}{2^*}\int_{\R^N}|u|^{2^*}dx,\quad u\in D^{1,2}(\R^N).
\ee
 We   consider the following modified problem to find  the nonnegative mountain pass solutions to (\ref{hhe6-4}),
\be\lab{hhe6-7}
\begin{cases}
-\Delta u-\mu_1 \frac{u}{|x|^2}=u_{+}^{2^*-1}+\alpha h(x)u_{+}^{\alpha-1}v_{+}^{\beta}, \; &\hbox{in}\;\R^N,\\
-\Delta v-\mu_2 \frac{v}{|x|^2}=v_{+}^{2^*-1}+\beta h(x)u_{+}^{\alpha}v_{+}^{\beta-1},\; &  \hbox{in}\;\R^N,
\end{cases}
\ee
where $u_+=\max\{u, 0\}$.  The  weak  solutions to problem (\ref{hhe6-7}) are critical points of the following functional $\overline{\Phi}:\D\rightarrow \R$ given by
\be\lab{zzzz-1}  \overline{\Phi}(u, v)=\overline{I}_{\mu_1}(u)+\overline{I}_{\mu_2}(v)-\int_{\R^N}h(x)u_{+}^{\alpha} v_{+}^{\beta} dx,\ee
where
\be\lab{hhe6-8}\overline{I}_{\mu_i}(w)=\frac{1}{2}Q_{\mu_i}(w)-\frac{1}{2^*}\int_{\R^N}w_{+}^{2^*}dx, \quad i=1,2,\ee
and
\be\lab{amy-e-1}
Q_\mu(w)=\int_{\R^N}|\nabla w|^2dx-\mu\int_{\R^N}\frac{w^2}{|x|^2}dx.
\ee
Obviously, the critical points of $\overline{\Phi}$ provide  nonnegative solutions to the original problem (\ref{hhe6-4}). We denote by $\overline{\mathcal{N}}$ the Nehari manifold associated to $\overline{\Phi}$, i.e.,
\be\lab{hhe6-8}
\overline{\mathcal{N}}=\Big\{(u, v)\in \D\backslash\{(0,0)\}:\langle \overline{\Phi}'(u, v), (u, v)\rangle=0\Big\}.
\ee
Assume $  N=3, \frac{1}{2}<\frac{1-4\mu_1}{1-4\mu_2}, \Theta \leq \min\{C_1, C_2\}$, where
$\Theta$ is defined in (\ref{ha1})  and
\be\lab{ling-e1}C_1:=\frac{(1-4\mu_1)^{\frac{2}{3}}-\big(\frac{1}{2}\big)^{\frac{2}{3}}\big(1-4\mu_2\big)^{\frac{2}{3}}}{ \big(\frac{1}{2}\big)^{\frac{\alpha-2}{6}}(1-4\mu_2)^{\frac{\alpha-2+4\beta}{6}}\alpha S^{\frac{\alpha-2}{4}+\beta-1}},\ee
\be\lab{ling-e2}C_2:=\frac{2-\sqrt[3]{2}}{ 2^{\frac{8-\beta}{6}}(1-4\mu_2)^{\frac{4\alpha+\beta-6}{6}}\beta S^{\frac{\beta-2}{4}+\alpha-1}}.\ee
It is easy to check that   $\displaystyle C_2>5\times 10^{-4}, C_1>10^{-3}\big[(1-4\mu_1)^{\frac{2}{3}}-\big(\frac{1}{2}\big)^{\frac{2}{3}}\big(1-4\mu_2\big)^{\frac{2}{3}}\big]$.

\bl\lab{hl-6-1}  Assume $N=3$, $\alpha, \beta\geq 2$,  $\alpha+\beta< 2^*$  and $\frac{1}{2}<\frac{1-4\mu_1}{1-4\mu_2}$. Let $\{(u_n,v_n)\}\subset\overline{\mathcal{N}}$ be a Palais-Smale sequence for $\overline{\Phi}|_{\overline{\mathcal{N}}}$ at level $c\in\R$. Then, there exists  $C>0$  such that  $||(u_n,v_n)||_{\D}\leq C$
for all $n\in \NN$ and $\overline{\Phi}'(u_n, v_n)\rightarrow 0$ in the dual space ${\D}^{*}$.  Furthermore, if   $\displaystyle\Theta\leq \min\{C_1, C_2\}$  and $c$ satisfies
\be\lab{hhe6-5}
\frac{1}{N}S^{\frac{N}{2}}(\mu_2)<c<\frac{1}{N}S^{\frac{N}{2}}(\mu_2)+\underset{(u, v)\in \overline{\mathcal{N}}}{\inf}\Phi(u, v);
\ee
\be\lab{hhe6-6}
c\neq \frac{l}{N}S^{\frac{N}{2}}(\mu_1)\; \hbox{ and }\; c\neq \frac{l}{N}S^{\frac{N}{2}}\;\hbox{for all }l\in \NN\backslash\{0\},
\ee
then up to a subsequence, $(u_n,v_n)\rightarrow (u_0, v_0)$ in $\D$.
\el.\br\lab{hl6-2}
In \cite[Lemma 3.5]{BVI}, the authors only considered the case  $\frac{1}{N}S^{\frac{N}{2}}(\mu_1)+\frac{1}{N}S^{\frac{N}{2}}(\mu_2)<\frac{1}{N}S^{\frac{N}{2}}$.  However, if $\mu_1+\mu_2\leq\frac{1}{4}$, then $\frac{1}{N}S^{\frac{N}{2}}(\mu_1)+\frac{1}{N}S^{\frac{N}{2}}(\mu_2)<\frac{1}{N}S^{\frac{N}{2}}$ will never meet. In particular,   the sign-changing $h(x)$ makes the proof  in  Lemma \ref{hl-6-1}  more complicated.
\er

\vskip 0.22in

\bp We divide the proof into five  steps.

\noindent {\bf Step 1:}
It is easy to show that  $\{(u_n, v_n)\}$ is bounded in $\D$ and $\overline{\Phi}'(u_n, v_n)\rightarrow 0$ in the dual space ${\D}^{*}$.  Up to a subsequence,   $\{(u_n ,v_n)\}_{n\in \NN}$  converges  weakly  to some $(u_0, v_0)$. Hence,   $((u_n)_-, (v_n)_-)\rightarrow (0,0)$ strongly in $\D$. It follows that $((u_n)_+, (v_n)_+)$ is a bounded Palais-Smale  sequence of $\overline{\Phi}$. For $((u_n)_+, (v_n)_+)$,  we may find a $t_n>0$ such that $t_n((u_n)_+, (v_n)_+)\in \mathcal{N}\cap \overline{\mathcal{N}}$. Since $(u_n, v_n)\in \overline{\mathcal{N}}$ and $((u_n)_-, (v_n)_-)\rightarrow 0$ in $\D$, we have $t_n\rightarrow 1$. Hence without loss of generality, we can assume  that
$\displaystyle u_n\geq 0, v_n\geq 0, \big\{(u_n, v_n)\big\}_{n\in\NN}\subset \mathcal{N}\cap \overline{\mathcal{N}}$ is a Palais-Smale sequence for $\overline{\Phi}$ at level $c$. Notice that $\Phi(u_n, v_n)=\overline{\Phi}(u_n, v_n)$.  For the simplicity, we use  $\displaystyle \|u\|_0$ to stand for $\|u\|_{D^{1,2}(\R^N)}$.
There exists $(u_0, v_0)\in {\D}$ and a subsequence, still denoted as $\displaystyle \big\{(u_n, v_n)\big\}_{n\in \NN}$ such that
\be\lab{lala1}
(u_n, v_n)\rightharpoonup (u_0, v_0)\quad \text{\rm{weakly in}} \; {\D},\quad (u_n, v_n)\rightarrow (u_0, v_0) \quad\text{a.e. in}\; \R^N,
\ee
\be\lab{lala2}
\text{and}\quad (u_n, v_n)\rightarrow (u_0, v_0) \quad \text{strongly in}\; L_{loc}^{\gamma}(\R^N)\times L_{loc}^{\gamma}(\R^N)\; \text{for all}\; \gamma\in [1, 2^{*}).
\ee
In view of the concentration-compactness principle due to Lions \cite{PL1, PL2}, there exists a subsequence, still denoted as $\big\{(u_n, v_n)\big\}_{n\in \NN}$, two at most countable sets $\mathcal{J}$ and $\mathcal{K}$,
set of points $\big\{x_j \in \R^N\backslash \{0\}:  j\in \mathcal{J}\big\}$ and $\big\{y_k \in \R^N\backslash \{0\}:  k\in \mathcal{K}\big\}$, real numbers $\zeta_j, \rho_j, j\in \mathcal{J}, \overline{\zeta}_k, \overline{\rho}_k, k\in \mathcal{K}, \zeta_0, \rho_0, \overline{\zeta}_0$ and $\overline{\rho}_0$ such that
\begin{equation}\lab{e6-11}
\begin{cases}
|\nabla u_n|^2 \rightharpoonup d\mu \geq |\nabla u_0|^2+\underset{j\in \mathcal{J}}{\sum}\zeta_j\delta_{x_j}+\zeta_0 \delta_0,\\
|\nabla v_n|^2 \rightharpoonup d\overline{\mu} \geq |\nabla v_0|^2+\underset{k \in \mathcal{K}}{\sum}\overline{\zeta}_k\delta_{y_k}+\overline{\zeta}_0 \delta_0,\\
|u_n|^{2^{*}} \rightharpoonup d\rho = |u_0|^{2^{*}}+\underset{j\in \mathcal{J}}{\sum}\rho_j\delta_{x_j}+\rho_0\delta_0,\\
|v_n|^{2^{*}} \rightharpoonup d\overline{\rho} = |v_0|^{2^{*}}+\underset{k\in \mathcal{K}}{\sum}\overline{\rho}_k\delta_{y_k}+\overline{\rho}_0\delta_0,\\
\frac{u_n^2}{|x|^2}\rightharpoonup d\theta=\frac{u_0^2}{|x|^2}+\theta_0\delta_0,\\
\frac{v_n^2}{|x|^2}\rightharpoonup d\overline{\theta}=\frac{v_0^2}{|x|^2}+\overline{\theta}_0\delta_0.
\end{cases}
\end{equation}
Define
\begin{align}\lab{eg}
&\zeta_\infty:=\lim_{R\rightarrow \infty}\limsup_{n\rightarrow \infty}\int_{|x|\geq R}|\nabla u_n|^2 dx,\;
\rho_\infty:=\lim_{R\rightarrow \infty} \limsup_{n \rightarrow \infty}\int_{|x|\geq R}|u_n|^{2^{*}}dx,\nonumber\\
&\overline{\zeta}_\infty:=\lim_{R\rightarrow \infty}\limsup_{n\rightarrow \infty}\int_{|x|\geq R}|\nabla v_n|^2 dx,\;
\overline{\rho}_\infty:=\lim_{R\rightarrow \infty} \limsup_{n \rightarrow \infty}\int_{|x|\geq R}|v_n|^{2^{*}}dx.
\end{align}
It follows that
\begin{align}\lab{eeg}
&\zeta_\infty:=\lim_{R\rightarrow \infty}\limsup_{n\rightarrow \infty}\int_{|x|\geq R}|\nabla (u_n-u_0)|^2 dx,\;\\
&\rho_\infty:=\lim_{R\rightarrow \infty} \limsup_{n \rightarrow \infty}\int_{|x|\geq R}|u_n-u_0|^{2^{*}}dx,\nonumber\\
&\overline{\zeta}_\infty:=\lim_{R\rightarrow \infty}\limsup_{n\rightarrow \infty}\int_{|x|\geq R}|\nabla (v_n-v_0)|^2 dx,\;\\
&\overline{\rho}_\infty:=\lim_{R\rightarrow \infty} \limsup_{n \rightarrow \infty}\int_{|x|\geq R}|v_n-v_0|^{2^{*}}dx.
\end{align}
From the Sobolev's inequality,  it follows easily that
\begin{equation}\lab{S1}
S\rho_{j}^{\frac{2}{2^{*}}}\leq \zeta_j \quad \text{for all}\; j\in \mathcal{J};\quad \quad
S\overline{\rho}_{k}^{\frac{2}{2^{*}}}\leq \overline{\zeta}_k\quad \text{for all}\; k \in \mathcal{K}.
\end{equation}

\vskip0.23in

 \noindent {\bf Step 2:}  We prove that either $u_n\rightarrow u_0$ strongly in $L^{2^*}(\R^N)$ or $v_n\rightarrow v_0$ strongly in $L^{2^*}(\R^N)$. If not, then there exist some $j_0\in\mathcal{J}\cup \{0,\infty\}$ and $k_0\in \mathcal{K}\cup\{0,\infty\}$ such that $\rho_{j_0}>0, \overline{\rho}_{k_0}>0$.  Since $\alpha+\beta<2^*$, we have $\rho_{_0}\geq S^{\frac{N}{2}}(\mu_1), \overline{\rho}_{k_0}\geq S^{\frac{N}{2}}(\mu_2)$.  In order to make the present paper easy to follow, we prefer to give part of the proofs.
 Indeed, for $\varepsilon>0$, let $\phi_j^\varepsilon$ be a smooth cut-off function centered at $x_j$, $0\leq \phi_j^\varepsilon \leq 1$, such that
\be\lab{t1}
\phi_j^\varepsilon(x)=
\begin{cases}
\hbox{1\quad if $|x-x_j|\leq \frac{\varepsilon}{2}$}\\
\hbox{0\quad if $|x-x_j|\geq \varepsilon$}
\end{cases}\hbox{and $\displaystyle |\nabla_y \phi_j^\varepsilon(y+x_j)|\leq \frac{4}{\varepsilon}$ for all $y=x-x_j\in \R^N$}.
\ee
Testing $\Phi'(u_n, v_n)$ with $(u_n \phi_j^\varepsilon, 0)$, we obtain
\begin{eqnarray}\lab{O}
0 &=& \lim_{n\rightarrow \infty}\big\langle\Phi^{'}(u_n, v_n), (u_n \phi_j^\varepsilon, 0)\big\rangle\nonumber\\
&=&\lim_{n\rightarrow \infty}\int_{\R^N}\Big(|\nabla u_n|^2\phi_j^\varepsilon+u_n\nabla u_n \cdot \nabla \phi_j^\varepsilon
-\phi_j^\varepsilon|u_n|^{2^{*}}  \nonumber\\
& & -\alpha h(x)|u_n|^\alpha |v_n|^\beta \phi_j^\varepsilon\Big)dx.
\end{eqnarray}
Notice that for all  $\varepsilon>0$ fixed,
\begin{eqnarray}\lab{O1}
\int_{\R^N}u_n\nabla u_n \cdot \nabla \phi_j^\varepsilon dx &=&\int_{\frac{\varepsilon}{2}\leq|x-x_j|\leq\varepsilon}u_n \nabla u_n \cdot \nabla \phi_j^\varepsilon dx\nonumber\\
&=&\int_{\frac{\varepsilon}{2}\leq|x-x_j|\leq\varepsilon}\big(u_n-u_0\big) \nabla u_n \cdot \nabla \phi_j^\varepsilon dx\nonumber\\
&&+\int_{\frac{\varepsilon}{2}\leq|x-x_j|\leq\varepsilon}u_0 \nabla u_n \cdot \nabla \phi_j^\varepsilon dx\nonumber\\
&:=&I+II.
\end{eqnarray}
Note $|\nabla \phi_j^\varepsilon|<\frac{4}{\varepsilon}$.  Without loss of generality,  we may assume that
\be\lab{lala2}
\quad (u_n, v_n)\rightarrow (u_0, v_0) \quad \text{strongly in}\; L_{loc}^{\gamma}(\R^N)\times L_{loc}^{\gamma}(\R^N)\; \text{for all}\; \gamma\in [1, 2^{*}).
\ee
Then by (\ref{lala2}) and the boundedness of $\nabla u_n$  in $L^2(\R^N)$, we have
\be\lab{lala3}
I\rightarrow 0 \quad \text{as}\; n\rightarrow +\infty.
\ee
Note   that  $II$ can be taken as a linear functional in $L^2(\R^N)$  and that  $u_n \rightharpoonup u_0$ in $D^{1,2}(\R^N)$, we have
\be\lab{lala4}
II\rightarrow \int_{\frac{\varepsilon}{2}\leq|x-x_j|\leq\varepsilon}u_0 \nabla u_0 \cdot \nabla \phi_j^\varepsilon dx \quad \text{as}\; n\rightarrow +\infty.
\ee
Since $(u_n, v_n)\rightharpoonup (u_0, v_0)$ weakly in ${\D}$, $(u_0, v_0)$ is a weakly solution to problem (\ref{hhe6-4}). Taking $(u_0\phi_j^\varepsilon, 0)$ as the testing function, then we have
\begin{eqnarray}\lab{lala5}
0=\int_{\R^N}\big[|\nabla u_0|^2\phi_j^\varepsilon+u_0\nabla u_0 \cdot \nabla \phi_j^\varepsilon
-\phi_j^\varepsilon|u_0|^{2^{*}}
-\alpha h(x)|u_0|^\alpha |v_0|^\beta \phi_j^\varepsilon\big]dx.
\end{eqnarray}
Further,
\begin{eqnarray*}
\int_{\R^N}|u_n|^\alpha |v_n|^\beta \phi_j^\varepsilon dx&=&\int_{|x-x_j|\leq \varepsilon}|u_n|^\alpha |v_n|^\beta \phi_j^\varepsilon dx\\
&\leq& \int_{|x-x_j|\leq \varepsilon}|u_n|^\alpha |v_n|^\beta dx\\
&\leq& \big(\int_{|x-x_j|\leq \varepsilon} |u_n|^{2^{*}}dx\big)^{\frac{\alpha}{2^{*}}}\big(\int_{|x-x_j|\leq \varepsilon} |v_n|^{2^{*}}dx\big)^{\frac{\beta}{2^{*}}}\\
&&\big(\int_{|x-x_j|\leq \varepsilon}1dx\big)^{\frac{2^{*}-\alpha-\beta}{2^{*}}}\\
&=& \hbox{$O(\varepsilon^{\frac{N(2^{*}-\alpha-\beta)}{2^{*}}})$.}
\end{eqnarray*}
Since $h(x)\in L^{\infty}(\R^N)$, we have that
\begin{equation}\lab{O2}
\int_{\R^N}\alpha h(x)|u_n|^\alpha |v_n|^\beta \phi_j^\varepsilon dx= O(\varepsilon^{\frac{N(2^{*}-\alpha-\beta)}{2^{*}}}).
\end{equation}
Especially,
\begin{equation}\lab{lala6}
\int_{\R^N}\alpha h(x)|u_0|^\alpha |v_0|^\beta \phi_j^\varepsilon dx= O(\varepsilon^{\frac{N(2^{*}-\alpha-\beta)}{2^{*}}}).
\end{equation}
By (\ref{e6-11}) and (\ref{O})$\sim$(\ref{lala6}), let $\varepsilon \rightarrow 0$ we obtain that
\be\lab{2013-8-1-e1}\mu_j - \rho_j \leq 0.\ee
By (\ref{S1}), we conclude that for all $j\in \mathcal{J}$, either $\rho_j=0$ or $\rho_j\geq S^{\frac{N}{2}}$, which also implies that $\mathcal{J}$ is finite. For the details about the similar results related to $j\in \{0,\infty\}$ and $\overline{\rho}_k, k\in \mathcal{K}\cup\{0,\infty\},$  we refer  the readers to \cite[Lemma 3.2]{ZZ}.  Then we have
    \begin{eqnarray*}
    c=\Phi(u_n, v_n)+o(1)&\geq& \big(\frac{1}{2}-\frac{1}{\alpha+\beta}\big)\big(S(\mu_1)\rho_{j_0}^{\frac{2}{2^*}}+S(\mu_2)\overline{\rho}_{k_0}^{\frac{2}{2^*}}\big)\\
    &&+\big(\frac{1}{\alpha+\beta}-\frac{1}{2^*}\big)(\rho_{j_0}+\overline{\rho}_{k_0})\\
    &\geq&\frac{1}{N}\big(S^{\frac{N}{2}}(\mu_1)+S^{\frac{N}{2}}(\mu_2)\big),
    \end{eqnarray*}
    which is a contradiction with (\ref{hhe6-5}).

 \vskip0.13in

 \noindent{\bf Step 3:}  We prove that either $u_n\rightarrow u_0$ strongly in $D^{1,2}(\R^N)$ or $v_n\rightarrow v_0$ strongly in $D^{1,2}(\R^N)$. By Step 2, if $u_n\rightarrow u_0$ strongly in $L^{2^*}(\R^N)$, then
    $$\|u_n-u_0\|_{\mu_1}^{2}=\langle \Phi'(u_n, v_n), (u_n-u_0, 0)\rangle+o(1)=o(1)\;\hbox{as}\;n\rightarrow +\infty.$$
    Hence, $u_n\rightarrow u_0$ strongly in $D^{1,2}(\R^N)$ in this case. If $v_n\rightarrow v_0$ strongly in $L^{2^*}(\R^N)$, correspondingly we have $v_n\rightarrow v_0$ strongly in $D^{1,2}(\R^N)$.

 \vskip0.13in

\noindent {\bf Step 4:} If $v_n\rightarrow v_0$ strongly in $D^{1,2}(\R^N)$, we prove    $u_n\to u_0$ strongly in $D^{1,2}(\R^N)$.
 We argue by contradiction and assume that $u_n\rightharpoonup u_0$ weakly but none of its subsequence  converges strongly to $u_0$.

Firstly we claim $v_0\not\equiv 0$. If not, $v_0\equiv 0$, it is easy to check that $\{u_n\}$ is a nonnegative Palais-Smale sequence for the functional  $I_{\mu_1}$  defined in (\ref{hhe6-2}), at the energy level $c=\underset{n\rightarrow\infty}{\lim} I_{\mu_1}(u_n)$,  which can be calculated as following:
\begin{align*}
c=&\underset{n\rightarrow \infty}{\lim} \Phi(u_n,v_n)\\
=&\underset{n\rightarrow \infty}{\lim} \frac{1}{2}\|(u_n, v_n)\|_{{\D}}^{2}-\frac{1}{2^{*}}\Big(|u_n|_{L^{2^{*}}(\R^N)}^{2^{*}}+|v_n|_{L^{2^{*}}(\R^N)}^{2^{*}}\Big)
-\int_{\R^N}h(x)|u_n|^{\alpha}|v_n|^{\beta}dx\\
=&\underset{n\rightarrow \infty}{\lim} I_{\mu_1}(u_n)+I_{\mu_2}(v_n)-\int_{\R^N}h(x)|u_n|^{\alpha}|v_n|^{\beta}dx\\
=&\underset{n\rightarrow \infty}{\lim} I_{\mu_1}(u_n)\;\hbox{since $v_n\rightarrow v_0\equiv 0$}.
\end{align*}
Then the result of \cite[Theorem 3.1]{SD} (we take $K(x)\equiv 1, \lambda=\mu_1$ in \cite{SD}) implies that there exists some $m, l\in \NN$ such that
$$c=\underset{n\rightarrow \infty}{\lim} \Phi(u_n,v_n)=\underset{n\rightarrow \infty}{\lim} I_{\mu_1}(u_n)=I_{\mu_1}(u_0)+\frac{m}{N}S^{\frac{N}{2}}+\frac{l}{N}S^{\frac{N}{2}}(\mu_1).$$
If $u_0\equiv 0$, then by (\ref{hhe6-6}), we have $m\neq 0, l\neq 0$.  In this case, $c\geq \frac{1}{N}S^{\frac{N}{2}}+\frac{1}{N}S^{\frac{N}{2}}(\mu_1)$, a contradiction with (\ref{hhe6-5}).
If $u_0\not\equiv 0$, since $u_0\geq 0$, we have $u_0=z_{\sigma}^{\mu_1}$ for some $\sigma>0$ and $\int_{\R^N}|u_0|^{2^*}dx=S^{\frac{N}{2}}(\mu_1), I_{\mu_1}(u_0)=\frac{1}{N}S^{\frac{N}{2}}(\mu_1)$. By (\ref{hhe6-6}), we obtain $m\neq 0$.  Then $c\geq \frac{1}{N}S^{\frac{N}{2}}+\frac{1}{N}S^{\frac{N}{2}}(\mu_1)$, also a contradiction with (\ref{hhe6-5}).
Thereby the claim $v_0\not\equiv 0$ is proved.

Thus we may assume that $u_n\rightharpoonup u_0$ weakly but not strongly in $D^{1,2}(\R^N)$ and $v_n\rightarrow v_0\not\equiv 0$ strongly in $D^{1,2}(\R^N)$.If $u_0\equiv 0$, then $v_0$ weakly solves
$$-\Delta v_0-\mu_2\frac{v_0}{|x|^2}-v_{0}^{2^*-1}=0,$$
 By  the known result in Section 1, we have $v_0=z_{\sigma}^{\mu_2}$ which is defined in (\ref{zzz-1}) for some $\sigma>0$. Thus, $|v_0|_{L^{2^*}(\R^N)}^{2^{*}}=\|v_0\|_{\mu_2}^{2}=S^{\frac{N}{2}}(\mu_2)$. Therefore,
\begin{align*}
c=&\big(\frac{1}{2}-\frac{1}{\alpha+\beta}\big)\|(u_n, v_n)\|_{\D}^2+\big(\frac{1}{\alpha+\beta}-\frac{1}{2^{*}}\big)\big(|u_n|_{L^{2^{*}}(\R^N)}^{2^{*}}+
|v_n|_{L^{2^{*}}(\R^N)}^{2^{*}}\big)+o(1)\\
\geq&\Big(\frac{1}{2}-\frac{1}{\alpha+\beta}\Big)\Big[\|(u_0, v_0)\|_{\D}^2 + \underset{j\in \mathcal{J}\cup\{0,\infty\}}{\Sigma}\zeta_j +\underset{k \in \mathcal{K}\cup\{0,\infty\} }{\Sigma}\overline{\zeta}_k\Big]+\\
&\Big(\frac{1}{\alpha+\beta}-\frac{1}{2^{*}}\Big)\Big[\int_{\R^N}(|u_0|^{2^{*}}+|v_0|^{2^{*}})dx +\underset{j\in \mathcal{J}\cup\{0,\infty\}}{\Sigma}\rho_j +\underset{k \in \mathcal{K}\cup\{0,\infty\} }{\Sigma}\overline{\rho}_k\Big]+o(1)\\
=&\Big(\frac{1}{2}-\frac{1}{\alpha+\beta}\Big)\|v_0\|_{\mu_2}^{2}+\Big(\frac{1}{\alpha+\beta}-\frac{1}{2^{*}}\Big)|v_0|_{L^{2^*}(\R^N)}^{2^{*}}\\
&+\Big(\frac{1}{\alpha+\beta}-\frac{1}{2^{*}}\Big)\big[\underset{j\in \mathcal{J}\cup \{0,\infty\}}{\Sigma}\zeta_j \big] +\Big(\frac{1}{\alpha+\beta}-\frac{1}{2^{*}}\Big)\big[\underset{j\in \mathcal{J}\cup\{0,\infty\}}{\Sigma}\rho_j\big]+o(1)\\
\geq&\frac{1}{N}|v_0|_{L^{2^*}(\R^N)}^{2^{*}}+\Big(\frac{1}{\alpha+\beta}-\frac{1}{2^{*}}\Big)\big[\underset{j\in \mathcal{J}}{\Sigma}S\rho_{j}^{\frac{2}{2^*}} +S(\mu_1)\rho_{0}^{\frac{2}{2^*}}+S(\mu_1)\rho_{\infty}^{\frac{2}{2^*}}\big]\\ &+\Big(\frac{1}{\alpha+\beta}-\frac{1}{2^{*}}\Big)\big[\underset{j\in \mathcal{J}}{\Sigma}\rho_j+\rho_0+\rho_\infty\big]+o(1).\\
\end{align*}
Since $u_n\rightharpoonup u_0$ weakly but not strongly in $D^{1,2}(\R^N)$, there exists some $j\in \mathcal{J}\cup \{0, \infty\}$ such that $\rho_j\neq 0$  and that
\begin{align*}
c\geq& \frac{1}{N}|v_0|_{L^{2^*}(\R^N)}^{2^{*}}+\Big(\frac{1}{\alpha+\beta}-\frac{1}{2^{*}}\Big)S(\mu_1)\rho_{j}^{\frac{2}{2^*}}
+\Big(\frac{1}{\alpha+\beta}-\frac{1}{2^{*}}\Big)\rho_j\\
\geq&\frac{1}{N}|v_0|_{L^{2^*}(\R^N)}^{2^{*}}+\Big(\frac{1}{\alpha+\beta}-\frac{1}{2^{*}}\Big)S(\mu_1)\cdot S^{\frac{N}{2}\cdot \frac{2}{2^*}}(\mu_1)+\Big(\frac{1}{\alpha+\beta}-\frac{1}{2^{*}}\Big)S^{\frac{N}{2}}(\mu_1)\\
=&\frac{1}{N}\big(S^{\frac{N}{2}}(\mu_1)+S^{\frac{N}{2}}(\mu_2)\big),
\end{align*}
a contradiction with (\ref{hhe6-5}). Hence we can assume that $u_0\not\equiv 0$. It is clear that $(u_0, v_0)\in \mathcal{N}\cap \overline{\mathcal{N}}$ and
\be\lab{531e-1}
\Phi(u_0, v_0)=\frac{1}{N}\Big(\int_{\R^N}u_{0}^{2^*} dx+\int_{\R^N}v_{0}^{2^*}dx\Big)+\frac{\alpha+\beta-2}{2}\int_{\R^N}h(x)u_0^\alpha v_0^\beta dx.
\ee
Since $v_n\rightarrow v_0$ in $D^{1,2}(\R^N)$ and $u_n\rightharpoonup u_0$ weakly in $D^{1,2}(\R^N)$, we can obtain that
\be\lab{531e-2}
\underset{n\rightarrow \infty}{\lim}\int_{\R^N} h(x)u_n^\alpha v_n^\beta dx=\int_{\R^N}h(x)u_0^\alpha v_0^\beta dx.
\ee
Combining  these facts, we have
\begin{align}\lab{61e-1}
\Phi(u_n, v_n)=&\frac{1}{2}\|(u_n, v_n)\|_{{\D}}^{2}-\frac{1}{2^{*}}\Big(|u_n|_{L^{2^{*}}(\R^N)}^{2^{*}}+|v_n|_{L^{2^{*}}(\R^N)}^{2^{*}}\Big)
-\int_{\R^N}h(x)|u_n|^{\alpha}|v_n|^{\beta}dx\nonumber\\
=&\frac{1}{N}\int_{\R^N}u_{n}^{2^*}dx+\frac{1}{N}\int_{\R^N}v_{n}^{2^*}dx+\frac{\alpha+\beta-2}{2}\int_{\R^N}h(x)u_n^\alpha v_n^\beta dx,
\end{align}
\be\lab{61e-2}
\underset{n\rightarrow \infty}{\lim}\int_{\R^N}v_{n}^{2^*}dx=\int_{\R^N}v_{0}^{2^*}dx,
\ee
\be\lab{61e-3}
\underset{n\rightarrow \infty}{\lim}\int_{\R^N}u_{n}^{2^*}dx=\int_{\R^N}u_{0}^{2^*}dx+\sum_{j\in \mathcal{J}}\rho_j+\rho_0+\rho_\infty.
\ee
By (\ref{531e-1}) $\sim$ (\ref{61e-3}), we obtain that
\begin{eqnarray}\lab{61e-4}
\Phi(u_0, v_0)&=&\underset{n\rightarrow \infty}{\lim}\Phi(u_n, v_n)-\frac{1}{N}\big(\sum_{j\in\mathcal{J}}\rho_j+\rho_0+\rho_\infty\big)\notag\\
&=&c-\frac{1}{N}\big(\sum_{j\in\mathcal{J}}\rho_j+\rho_0+\rho_\infty\big)\notag\\
&\leq &c-\frac{1}{N}S^{\frac{N}{2}}(\mu_1)\notag\\
&<&\frac{1}{N}S^{\frac{N}{2}}(\mu_2).
\end{eqnarray}
Note  that
\be\lab{amy-e-5}\int_{\R^N}u_{0}^{2^*}dx+\alpha\int_{\R^N}h(x)u_0^\alpha v_0^\beta dx=\|u_0\|_{\mu_1}^{2}\geq S(\mu_1)\big(\int_{\R^N}u_{0}^{2^*}dx\big)^{\frac{2}{2^*}}. \ee
If $\displaystyle\int_{\R^N}h(x)u_0^\alpha v_0^\beta dx\leq 0$, then $\displaystyle\int_{\R^N}u_{0}^{2^*}dx\geq S^{\frac{N}{2}}(\mu_1)$.  Similarly, we have $\displaystyle\int_{\R^N}v_{0}^{2^*}dx\geq S^{\frac{N}{2}}(\mu_2)$.
Thus,
\begin{eqnarray}\lab{61e-5}
\Phi(u_0, v_0)&=&(\frac{1}{2}-\frac{1}{\alpha+\beta})\|(u_0, v_0)\|_{\D}^{2}+(\frac{1}{\alpha+\beta}-\frac{1}{2^*})\big(\|u_0\|_{L^{2^*}}^{2^*}+\|v_0\|_{L^{2^*}}^{2^*}\big)\notag\\
&\geq&(\frac{1}{2}-\frac{1}{\alpha+\beta})S(\mu_1)\|u_0\|_{L^{2^*}}^{2}+(\frac{1}{2}-\frac{1}{\alpha+\beta})S(\mu_2)\|v_0\|_{L^{2^*}}^{2}\notag\\
&&+(\frac{1}{\alpha+\beta}-\frac{1}{2^*})\big(\|u_0\|_{L^{2^*}}^{2^*}+\|v_0\|_{L^{2^*}}^{2^*}\big)\notag\\
&\geq&(\frac{1}{2}-\frac{1}{\alpha+\beta})S(\mu_1)S^{\frac{N}{2^*}}(\mu_1)+(\frac{1}{\alpha+\beta}-\frac{1}{2^*})S^{\frac{N}{2}}(\mu_1)\notag\\
&&+(\frac{1}{2}-\frac{1}{\alpha+\beta})S(\mu_2)S^{\frac{N}{2^*}}(\mu_2)+(\frac{1}{\alpha+\beta}-\frac{1}{2^*})S^{\frac{N}{2}}(\mu_2)\notag\\
&=&\frac{1}{N}\big(S^{\frac{N}{2}}(\mu_1)+S^{\frac{N}{2}}(\mu_2)\big)>\frac{1}{N}S^{\frac{N}{2}}(\mu_2),
\end{eqnarray}
a contradiction with (\ref{61e-4}).

If $\displaystyle\int_{\R^N}h(x)u_0^\alpha v_0^\beta dx>0$, then by (\ref{531e-1}) and (\ref{61e-4}) we deduce that
\be\lab{61e-6}
\|u_0\|_{L^{2^*}(\R^N)}^{2^*}+\|v_0\|_{L^{2^*}(\R^N)}^{2^*}<S^{\frac{N}{2}}(\mu_2).
\ee
Similar to (\ref{amy-e-5}), we have
\be\lab{amy-e-6}\int_{\R^N}v_{0}^{2^*}dx+\beta\int_{\R^N}h(x)u_0^\alpha v_0^\beta dx=\|v_0\|_{\mu_2}^{2}\geq S(\mu_2)\big(\int_{\R^N}v_{0}^{2^*}dx\big)^{\frac{2}{2^*}}. \ee
Then by (\ref{61e-6}), (\ref{amy-e-5}) and (\ref{amy-e-6}) we obtain that
\be\lab{amy-e-4}
S(\mu_2)\|v_0\|_{L^{2^*}(\R^N)}^{2}\leq\|v_0\|_{L^{2^*}(\R^N)}^{2^*}+ \beta \Theta S^{\frac{N-2}{4}\alpha}(\mu_2)\|v_0\|_{L^{2^*}(\R^N)}^{\beta}
\ee
and
\be\lab{amy-e-7}
S(\mu_1)\|u_0\|_{L^{2^*}(\R^N)}^{2}\leq\|u_0\|_{L^{2^*}(\R^N)}^{2^*}+ \alpha \Theta S^{\frac{N-2}{4}\beta}(\mu_2)\|u_0\|_{L^{2^*}(\R^N)}^{\alpha}.
\ee
Since $N=3$, we have
$ S(\mu)=(1-4\mu)^{\frac{2}{3}}S,$
then (\ref{amy-e-4}) and (\ref{amy-e-7})  are  equivalent to
\be\lab{amy-e-8}
\big(1-4\mu_2\big)^{\frac{2}{3}}S\leq \|v_0\|_{L^6(\R^3)}^{4}+\beta \Theta\big(1-4\mu_2\big)^{\frac{2\alpha}{3}}S^\alpha\|v_0\|_{L^6(\R^3)}^{\beta-2}
\ee
and
\be\lab{amy-e-9}
\big(1-4\mu_1\big)^{\frac{2}{3}}S\leq \|u_0\|_{L^6(\R^3)}^{4}+\alpha \Theta\big(1-4\mu_2\big)^{\frac{2\beta}{3}}S^\beta\|u_0\|_{L^6(\R^3)}^{\alpha-2}.
\ee
Note   that  $\alpha\geq 2$ and  that  $f(t):=t^{\frac{2}{3}}+\alpha\Theta(1-4\mu_2)^{\frac{2}{3}\beta}S^\beta t^{\frac{\alpha-2}{6}}$ is increasing  in $(0,+\infty)$. If $\Theta\leq C_1$, where $C_1$ is given in (\ref{ling-e1}), then
\begin{eqnarray*}
&&f\big(\frac{1}{2}(1-4\mu_2)S^{\frac{3}{2}}\big)\\
& &=\big[\frac{1}{2}(1-4\mu_2)S^{\frac{3}{2}}\big]^{\frac{2}{3}}+\alpha\Theta(1-4\mu_2)^{\frac{2}{3}\beta}S^\beta
\big[\frac{1}{2}(1-4\mu_2)S^{\frac{3}{2}}\big]^{\frac{\alpha-2}{6}}\\
& &=\big(\frac{1}{2}\big)^{\frac{2}{3}}\big(1-4\mu_2\big)^{\frac{2}{3}}S+\big(\frac{1}{2}\big)^{\frac{\alpha-2}{6}}
\alpha\Theta\big(1-4\mu_2\big)^{\frac{2}{3}\beta}\big(1-4\mu_2\big)^{\frac{\alpha-2}{6}}S^{\beta+\frac{\alpha-2}{4}}\\
& & \leq  \big(1-4\mu_1\big)^{\frac{2}{3}}S.
\end{eqnarray*}
Recalling (\ref{amy-e-9}),  it follows that
$ \|u_0\|_{L^{2^*}(\R^N)}^{2^*}\geq \frac{1}{2}(1-4\mu_2)S^{\frac{3}{2}}=\frac{1}{2}S^{\frac{N}{2}}(\mu_2).$
Similarly,   we have
$ \|v_0\|_{L^{2^*}(\R^N)}^{2^*}\geq \frac{1}{2}S^{\frac{N}{2}}(\mu_2).$
Thus,
$ \|u_0\|_{L^{2^*}(\R^N)}^{2^*}+\|v_0\|_{L^{2^*}(\R^N)}^{2^*}\geq S^{\frac{N}{2}}(\mu_2),$
a contradiction to (\ref{61e-6}).

\vskip0.13in
\noindent{\bf Step 5:}   If $u_n\rightarrow u_0$ strongly in $D^{1,2}(\R^N)$, we show that  $v_n$ converges strongly. For  this case,  $\mathcal{K}\cup \{0, \infty\} $ is reduced to be at most one point. In fact,  if not, by (3.31), (3.39), we obtain that
$ c\geq \frac{2}{N}S^{\frac{N}{2}}(\mu_2)\geq\frac{1}{N}\big(S^{\frac{N}{2}}(\mu_1)+S^{\frac{N}{2}}(\mu_2)\big),$
  contradicting to (\ref{hhe6-5}). Assume that $v_n\rightharpoonup v_0$ weakly but none of its subsequence  converges strongly to $v_0$. We claim $u_0\not\equiv 0$. If not, $u_0\equiv 0$ and $v_0\not\equiv 0$, then $v_0\in D^{1,2}(\R^N)$ is a weak solution to
    \be\lab{he6-3}-\Delta v_0-\mu_2\frac{v_0}{|x|^2}=|v_0|^{2^*-1}v_0.\ee
    Denote
    $$E_1:=\inf\big\{ I_{\mu_2}(v)\big|v\neq 0\;\hbox{is a solution of problem (\ref{he6-3})}\big\},$$
    $$E_2:=\inf\big\{ \int_{\R^N}|v|^{2^*}dx\big|v\neq 0\;\hbox{is a solution of problem (\ref{he6-3})}\big\},$$ where $I_\mu(v)$ is defined in (\ref{hhe6-2}).
    Then it is well known that $E_1=\frac{1}{N}S^{\frac{N}{2}}(\mu_2), E_2=S^{\frac{N}{2}}(\mu_2)$, and they are only achieved by
      $$
     \begin{cases}
     &\pm z_{\sigma}^{\mu_2}, \sigma>0\;\hbox{ if}\; \mu_2>0\\
     & \pm z_{\sigma, x_i}, \sigma>0, x_i\in \R^N \;\hbox{if}\;\mu_2=0.
     \end{cases}
     $$
      Thus, by the fact  that  $\{v_n\}$ concentrate at exactly one point, we obtain that
    $$c\geq \frac{1}{N}\Big(\int_{\R^N}v_{0}^{2^*}dx+S^{\frac{N}{2}}(\mu_2)\Big)\geq\frac{2}{N}S^{\frac{N}{2}}(\mu_2)\geq \frac{1}{N}(S^{\frac{N}{2}}(\mu_1)+S^{\frac{N}{2}}(\mu_2)),$$
    contradicting (\ref{hhe6-5}).
 On the other hand, if $u_0\equiv 0, v_0\equiv 0$, then $v_n$ solves
    \be\lab{530e-1}
    -\Delta v_n-\mu_2\frac{v_n}{|x|^2}-v_{n}^{2^*-1}=o(1),\quad\hbox{in the dual space} \big(D^{1,2}(\R^N)\big)^*.
    \ee
   Then an analogous argument as that in step 4  will lead to that
   \begin{align*}
   c&=\Phi(u_n, v_n)+o(1)=I_{\mu_2}(v_n)+o(1)\rightarrow I_{\mu_2}(v_0)+\frac{m}{N}S^{\frac{N}{2}}(\mu_2)+ \frac{l}{N}S^{\frac{N}{2}}\\
   &=\frac{m}{N}S^{\frac{N}{2}}(\mu_2)+\frac{l}{N}S^{\frac{N}{2}}\;\hbox{for some $m, l\in \NN$}
   \end{align*}
as $n\rightarrow \infty$. By (\ref{hhe6-6}), we have $m\neq 0, l\neq 0$, then $c\geq \frac{1}{N}S^{\frac{N}{2}}(\mu_2)+\frac{1}{N}S^{\frac{N}{2}}$, a contradiction with (\ref{hhe6-5}).
Thereby the claim $u_0\not\equiv 0$ is proved. Thus we may assume that $u_n\rightarrow u_0\not\equiv 0$ strongly in $D^{1,2}(\R^N)$, and $v_n\rightharpoonup v_0$.  If $v_0\equiv 0$, then $u_0$ weakly solves
 $$-\Delta u_0-\mu_1 \frac{u_0}{|x|^2}-u_{0}^{2^*-1}=0,$$
 hence $u_0=z_{\sigma}^{\mu_1}$ for some $\sigma>0$ and $\displaystyle \int_{\R^N} |u_0|^{2^*}dx=S^{\frac{N}{2}}(\mu_1)$. When $\mu_2>0$ (similarly for the case of $\mu_2=0$), by the fact that  $\{v_n\}$ concentrates  to exactly one point, we deduce that
 \begin{eqnarray}\lab{e6-25}
c&=&\Big(\frac{1}{2}-\frac{1}{\alpha+\beta}\Big)\Big[\|u_0\|_{\mu_1}^{2}+\|v_0\|_{\mu_2}^2 \nonumber\\
&&+\underset{k \in \mathcal{K} }{\sum}\overline{\zeta}_k+(\overline{\zeta}_0-\mu_2\overline{\theta}_0)+(\overline{\zeta}_\infty-\mu_2\overline{\theta}_\infty)\Big]+\Big(\frac{1}{\alpha+\beta}-\frac{1}{2^{*}}\Big)\nonumber\\
&&\Big[\int_{\R^N}(|u_0|^{2^{*}}+|v_0|^{2^{*}})dx +\underset{k \in \mathcal{K}}{\sum}\overline{\rho}_k+\overline{\rho}_0+\overline{\rho}_\infty\Big]+o(1)\nonumber\\
&\geq & \Big(\frac{1}{2}-\frac{1}{\alpha+\beta}\Big)\Big[\|u_0\|_{\mu_1}^{2}+S\underset{k \in \mathcal{K}}{\sum}\overline{\rho}_{k}^{\frac{2}{2^{*}}}
+S(\mu_2)(\overline{\rho}_{0}^{\frac{2}{2^*}}+\overline{\rho}_{\infty}^{\frac{2}{2^*}})\Big]\nonumber\\
&&+\Big(\frac{1}{\alpha+\beta}-\frac{1}{2^{*}}\Big)\Big(\|u_0\|_{2^*}^{2^*}+\underset{k \in \mathcal{K}}{\sum}\overline{\rho}_k+\overline{\rho}_0+\overline{\rho}_\infty\Big)+o(1)\nonumber\\
&\geq &\Big(\frac{1}{2}-\frac{1}{\alpha+\beta}\Big)\Big(S(\mu_1)\|u_0\|_{2^*}^{2}+S(\mu_2)\overline{\rho}_{\widetilde{k}}^{\frac{2}{2^*}}\Big)\nonumber\\
&&+\Big(\frac{1}{\alpha+\beta}-\frac{1}{2^{*}}\Big)\Big(\|u_0\|_{2^*}^{2^*}+\overline{\rho}_{\widetilde{k}}\Big)\quad \hbox{for some $\widetilde{k}\in \mathcal{K}\cup \{0,\infty\}$}+o(1)\nonumber\\
&\geq& \frac{1}{N}\big(S^{\frac{N}{2}}(\mu_1)+S^{\frac{N}{2}}(\mu_2)\big)+o(1),
\end{eqnarray}
 a contradiction with (\ref{hhe6-5}). Hence,  both $u_0\not\equiv 0$ and $v_0\not\equiv 0$.  Hence,
 \begin{eqnarray}
 c&=&\Phi(u_n, v_n)-\frac{1}{2}\big\langle\Phi'(u_n, v_n), (u_n, v_n)\big\rangle+o(1)\\
 &=&\frac{1}{N}\Big(\int_{\R^N} u_{n}^{2^*}dx+\int_{\R^N}v_{n}^{2^*}dx\Big)+\frac{\alpha+\beta-2}{2}\int_{\R^N}h(x)u_n^\alpha v_n^\beta dx+o(1),\notag
 \end{eqnarray}
it follows that, for some $k\in \mathcal{K}\cup \{0,\infty\}$,
 \be\lab{530e-2}
 c=\frac{1}{N}\big(\int_{\R^N}u_{0}^{2^*} dx+\int_{\R^N}v_{0}^{2^*}dx+\overline{\rho}_k\big)+\frac{\alpha+\beta-2}{2}\int_{\R^N} h(x)u_0^\alpha v_0^\beta dx.
 \ee
 On the other hand,  recall that   $\langle \Phi'(u_n, v_n), (u_0, v_0)\rangle=o(1)$, we have
 $$\|(u_0, v_0)\|_{\D}^{2}=\int_{\R^N}u_{0}^{2^*}dx+\int_{\R^N}v_{0}^{2^*}dx+(\alpha+\beta)\int_{\R^N} h(x)u_0^\alpha v_0^\beta dx, $$
 i.e., $(u_0, v_0)\in \mathcal{N}$. If $\mu_2>0$, by (3.31),  (\ref{530e-2})  and assumption (\ref{hhe6-5}), we obtain that
 \begin{eqnarray}
 \Phi(u_0, v_0)&=&\frac{1}{N}\Big(\int_{\R^N} u_{0}^{2^*} dx +\int_{\R^N}v_{0}^{2^*}dx\Big)+\frac{\alpha+\beta-2}{2}\int_{\R^N}h(x)u_0^\alpha v_0^\beta dx\notag\\
 &=& c-\frac{\overline{\rho}_k}{N}<\frac{1}{N}S^{\frac{N}{2}}(\mu_2)+\underset{(u,v)\in\mathcal{N}}{\inf}\Phi(u, v)-\frac{\overline{\rho}_k}{N}\notag\\
 &\leq &\frac{1}{N}S^{\frac{N}{2}}(\mu_2)+\underset{(u,v)\in\mathcal{N}}{\inf}\Phi(u, v)-\frac{1}{N}S^{\frac{N}{2}}(\mu_2)\notag\\
 &=&\underset{(u,v)\in\mathcal{N}}{\inf}\Phi(u, v),
 \end{eqnarray}
  a contradiction with $(u_0, v_0)\in \mathcal{N}$.
 \ep

\vskip0.336in
Next, we consider the case of $\alpha+\beta=2^*$.  For this case,   we can not expect $\rho_j=0$ or $\rho_j\geq S^{\frac{N}{2}}$ for $j\not\in\{0,\infty\}$  any more. Because the H\"{o}lder's inequality is not enough to ensure
\be\lab{ml-e1}
\int_{\R^N}h(x)|u_n|^\alpha |v_n|^\beta \phi_j^\varepsilon dx \rightarrow 0\;\hbox{uniformly as}\; \varepsilon \rightarrow 0,
\ee
where $\phi_j^\varepsilon$ is defined in (\ref{t1}).  Thus, the step 1 of Lemma \ref{hl-6-1}   fails  for the case of $\alpha+\beta=2^*$. We will impose  more conditions on $h(x)$ to overcome this difficulty.  Assume $\mu_1+\mu_2\neq 0$, take $\varepsilon_1$ small enough such that
\be\lab{ml-e11}
2(1-\varepsilon_1)^{\frac{N-2}{N}}>\big(1-\frac{4\mu_1}{(N-2)^2}\big)^{\frac{N-1}{2}}+\big(1-\frac{4\mu_2}{(N-2)^2}\big)^{\frac{N-1}{2}}.
\ee
For example, if $N=3$, we  choose  $\varepsilon_1$ satisfying
\be\lab{ml-e12}
\varepsilon_1<1-\big[1-2(\mu_1+\mu_2)\big]^3
\ee
and if $N=4$, we may take $\varepsilon_1$ satisfying
\be\lab{ml-e13}
\varepsilon_1<1-\frac{\big[(1-\mu_1)^{\frac{3}{2}}+(1-\mu_2)^{\frac{3}{2}}\big]^2}{4}.
\ee

\bl\lab{lemma-ml1}
Assume $\alpha+\beta= 2^*$, $0\leq \mu_2\leq \mu_1<\frac{1}{4}, 1-4\mu_2<2-8\mu_1, \alpha, \beta\geq 2,  \mu_1+\mu_2\neq 0$. Let $\{(u_n,v_n)\}\subset\overline{\mathcal{N}}$ be a Palais-Smale sequence for $\overline{\Phi}|_{\overline{\mathcal{N}}}$ at level $c\in\R$. Then, there exists some constant $C$  such that  $||(u_n,v_n)||_{\D}\leq C$
for all $n\in \NN$ and $\overline{\Phi}'(u_n, v_n)\rightarrow 0$ in the dual space ${\D}^{*}$. Moreover, if $c$ satisfies (\ref{hhe6-5}), (\ref{hhe6-6}) and
\be\lab{mla-e1}
\Theta\leq \min\Big\{\frac{1-(1-\varepsilon_1)^{\frac{2}{N}}}{2^{\frac{\beta}{2}}\alpha (1-\varepsilon_1)^{\frac{\alpha-2}{2^*}}}, \frac{1-(1-\varepsilon_1)^{\frac{2}{N}}}{2^{\frac{\alpha}{2}}\beta (1-\varepsilon_1)^{\frac{\beta-2}{2^*}}}, C_1, C_2\Big\},
\ee where $\varepsilon_1$ satisfies (\ref{ml-e12}) and $C_1, C_2$  are defined in  (\ref{ling-e1}) and (\ref{ling-e2}),
then, up to a subsequence, $(u_n,v_n)\rightarrow (u_0, v_0)$ in $\D$.
\el
\br\lab{lj-r1}
Under the assumptions of Lemma \ref{lemma-ml1}, we can have
\begin{align*}
&\min\Big\{\frac{1-(1-\varepsilon_1)^{\frac{2}{N}}}{2^{\frac{\beta}{2}}\alpha (1-\varepsilon_1)^{\frac{\alpha-2}{2^*}}}, \frac{1-(1-\varepsilon_1)^{\frac{2}{N}}}{2^{\frac{\alpha}{2}}\beta (1-\varepsilon_1)^{\frac{\beta-2}{2^*}}}, C_1, C_2\Big\}>\\
&\min\Big\{\frac{\mu_1+\mu_2-(\mu_1+\mu_2)^2}{6}, 10^{-3}\big[(1-4\mu_1)^{\frac{2}{3}}-\big(\frac{1}{2}\big)^{\frac{2}{3}}\big(1-4\mu_2\big)^{\frac{2}{3}}\big],5\times 10^{-4}\Big\}.
\end{align*}
\er
\bp We need several steps.

\noindent{\bf Step 1:}
There exist  an at  most countable set  $\mathcal{J}$ (for simplicity, here we view $\mathcal{J}$ as the set  $\mathcal{J}\cup \mathcal{K}$ in Lemma \ref{hl-6-1}), the set of points $\big\{x_j \in \R^N\backslash \{0\}:  j\in J\big\}$, real numbers $\zeta_j, \rho_j, \overline{\zeta}_j, \overline{\rho}_j, j\in J, \zeta_0, \rho_0, \overline{\zeta}_0$ and $\overline{\rho}_0$,  such that
\begin{equation}\lab{ml-e5}
\begin{cases}
|\nabla u_n|^2 \rightharpoonup d\mu \geq |\nabla u_0|^2+\underset{j\in \mathcal{J}}{\sum}\zeta_j\delta_{x_j}+\zeta_0 \delta_0,\\
|\nabla v_n|^2 \rightharpoonup d\overline{\mu} \geq |\nabla v_0|^2+\underset{j \in \mathcal{J}}{\sum}\overline{\zeta}_j\delta_{x_j}+\overline{\zeta}_0 \delta_0,\\
|u_n|^{2^{*}} \rightharpoonup d\rho = |u_0|^{2^{*}}+\underset{j\in \mathcal{J}}{\sum}\rho_j\delta_{x_j}+\rho_0\delta_0,\\
|v_n|^{2^{*}} \rightharpoonup d\overline{\rho} = |v_0|^{2^{*}}+\underset{j\in \mathcal{J}}{\sum}\overline{\rho}_j\delta_{x_j}+\overline{\rho}_0\delta_0,\\
\frac{u_n^2}{|x|^2}\rightharpoonup d\theta=\frac{u_0^2}{|x|^2}+\theta_0\delta_0,\\
\frac{v_n^2}{|x|^2}\rightharpoonup d\overline{\theta}=\frac{v_0^2}{|x|^2}+\overline{\theta}_0\delta_0.
\end{cases}
\end{equation}
Note $\Phi'(u_n, v_n)\rightarrow 0$ in $D^*$ and $\{(u_n, v_n)\}$ is bounded,  we obtain that
\begin{eqnarray}\lab{ml-e4}
0&=&\underset{n\rightarrow \infty}{\lim}\big\langle \Phi'(u_n, v_n), (u_n\phi_j^\varepsilon, 0) \big\rangle\nonumber\\
&=&\underset{n\rightarrow \infty}{\lim}\int_{\R^N}\Big(|\nabla u_n|^2\phi_j^\varepsilon+u_n\nabla u_n \cdot \nabla \phi_j^\varepsilon-\mu_1\frac{u_n^2\phi_j^\varepsilon}{|x|^2}\nonumber\\
&&-\phi_j^\varepsilon |u_n|^{2^*}-\alpha h(x)|u_n|^\alpha |v_n|^\beta \phi_j^\varepsilon\Big)dx.
\end{eqnarray}
 By $(H'_1)$, we can follow  the process of \cite[Lemma 3.2]{ZZ} and obtain that
\be\lab{ml-e2}
\rho_0=0\;\hbox{or}\;\rho_0\geq S^{\frac{N}{2}}(\mu_1),\;\;  \rho_\infty=0\;\hbox{or}\;\rho_\infty\geq S^{\frac{N}{2}}(\mu_1).
\ee
Similarly we can obtain that
\be\lab{ml-e3}
\overline{\rho}_0=0\;\hbox{or}\;\overline{\rho}_0\geq S^{\frac{N}{2}}(\mu_2),\;\; \overline{\rho}_\infty=0\;\hbox{or}\;\overline{\rho}_\infty\geq S^{\frac{N}{2}}(\mu_2).
\ee
For $x_j\in \R^N\backslash\{0\}$, if $h(x_j)\leq 0$, then we can argue as above and obtain that
\be\lab{ml-e7}
\rho_j=0\;\hbox{or}\;\rho_j\geq S^{\frac{N}{2}}; \overline{\rho}_j=0\;\hbox{or}\;\overline{\rho}_j\geq S^{\frac{N}{2}}.
\ee
Next we consider $x_j\in \R^N\backslash\{0\}$ with $h(x_j)>0$.
Then one of the following holds:
\begin{itemize}
\item[(1)] $\rho_j=\overline{\rho}_j=0$;
\item[(2)] $\rho_j=0$ and $\overline{\rho}_j>0$;
\item[(3)] $\rho_j>0$ and $\overline{\rho}_j=0$;
\item[(4)] $\rho_j>0$ and $\overline{\rho}_j>0$.
\end{itemize}
If (3) holds, then (\ref{ml-e1}) satisfies.  By (\ref{ml-e4}) and $x_j\neq 0$, we can obtain
$\zeta_j-\rho_j\leq 0.$
Then the Sobolev's  inequality implies that either  $ \rho_j=0\;\hbox{or}\;\rho_j\geq S^{\frac{N}{2}}.$
Thus, if (3) holds, we have $\rho_j\geq S^{\frac{N}{2}}$. Similarly, if (2) holds, we have $\overline{\rho}_j\geq S^{\frac{N}{2}}$.
If (4) holds, we can not apply (\ref{ml-e1}) but we claim
 \be\lab{mla-e2}
S(\rho_{j}^{\frac{2}{2^*}}+\overline{\rho}_{j}^{\frac{2}{2^*}})\geq S^{\frac{N}{2}}(\mu_1)+S^{\frac{N}{2}}(\mu_2).
\ee
Without loss of generality, we can assume that $\max\{\rho_j, \overline{\rho}_j\}\leq 2^{\frac{2^*}{2}}S^{\frac{N}{2}}$.   Then by the  Sobolev's inequality and (\ref{ml-e4}),  we have
\be\lab{ml-e9}
S \rho_{j}^{\frac{2}{2^*}}-\rho_j-\alpha \Theta \rho_{j}^{\frac{\alpha}{2^*}}(2^{\frac{2^*}{2}}S^{\frac{N}{2}})^{\frac{\beta}{2^*}}\leq 0,
\ee
which  is equivalent to
$S\leq \rho_{j}^{1-\frac{2}{2^*}}+\alpha \Theta 2^{\frac{\beta}{2}} S^{\frac{(N-2)\beta}{4}} \rho_{j}^{\frac{\alpha-2}{2^*}}.$
Consider the function $\displaystyle f(t)=t^{1-\frac{2}{2^*}}+\alpha\Theta 2^{\frac{\beta}{2}}S^{\frac{(N-2)\beta}{4}}t^{\frac{\alpha-2}{2^*}}$, which is increasing in $(0, +\infty)$ because of $\alpha\geq 2$.
If $\displaystyle \Theta\leq \frac{1-(1-\varepsilon_1)^{\frac{2}{N}}}{2^{\frac{\beta}{2}}\alpha (1-\varepsilon_1)^{\frac{\alpha-2}{2^*}}}$, then we can compute $f\big((1-\varepsilon_1)S^{\frac{N}{2}}\big)$ as following:
\begin{eqnarray}\lab{ml-e11}
f\big((1-\varepsilon_1)S^{\frac{N}{2}}\big)
 = \big[(1-\varepsilon_1)^{\frac{2}{N}}+\alpha\Theta2^{\frac{\beta}{2}}(1-\varepsilon_1)^{\frac{\alpha-2}{2^*}}\big]S
 \leq S,  &
\end{eqnarray}
which implies that $\rho_j\geq (1-\varepsilon_1)S^{\frac{N}{2}}$ and $\displaystyle S \rho_{j}^{\frac{2}{2^*}}\geq (1-\varepsilon_1)^{\frac{N-2}{N}}S^{\frac{N}{2}}$.
Similarly,  if $\displaystyle \Theta\leq \frac{1-(1-\varepsilon_1)^{\frac{2}{N}}}{2^{\frac{\alpha}{2}}\beta (1-\varepsilon_1)^{\frac{\beta-2}{2^*}}}$, we have $\displaystyle S \overline{\rho}_{j}^{\frac{2}{2^*}}\geq (1-\varepsilon_1)^{\frac{N-2}{N}}S^{\frac{N}{2}}$. Hence,  by (\ref{ml-e11}) we have
$$S(\rho_{j}^{\frac{2}{2^*}}+\overline{\rho}_{j}^{\frac{2}{2^*}})\geq 2(1-\varepsilon_1)^{\frac{N-2}{N}}S^{\frac{N}{2}}\geq S^{\frac{N}{2}}(\mu_1)+S^{\frac{N}{2}}(\mu_2).$$
\vskip 0.130in
 \noindent{\bf Step 2: } We prove that either $u_n\rightarrow u_0$ or  $v_n\rightarrow v_0$ strongly in $L^{2^*}(\R^N)$. If not, then  there exist some $j_0\in\mathcal{J}\cup \{0\}\cup \{\infty\}$ and $j_1\in \mathcal{J}\cup\{0\}\cup\{\infty\}$ such that $\rho_{j_0}>0, \overline{\rho}_{j_1}>0$.   If $j_0\neq j_1$, then we have $\rho_{j_0}\geq S^{\frac{N}{2}}(\mu_1), \overline{\rho}_{j_1}\geq S^{\frac{N}{2}}(\mu_2)$ and obtain that
    \begin{eqnarray*}
    c&=&\Phi(u_n, v_n)+o(1)\\
    &=&\frac{1}{N}\|(u_n, v_n)\|_{\D}^{2}+o(1)\\
    &\geq& \frac{1}{N} S(\mu_1)\rho_{j_0}^{\frac{2}{2^*}}+\frac{1}{N}S(\mu_2)\overline{\rho}_{j_1}^{\frac{2}{2^*}}\\
    &\geq& \frac{1}{N}\big(S^{\frac{N}{2}}(\mu_1)+S^{\frac{N}{2}}(\mu_2)\big),
    \end{eqnarray*}
    which   contradicts  with (\ref{hhe6-5}).
    If $j_0=j_1=j$  and $j\in\{0,\infty\}$, then we have $\rho_{j}\geq S^{\frac{N}{2}}(\mu_1), \overline{\rho}_{j}\geq S^{\frac{N}{2}}(\mu_2)$ and $\displaystyle c\geq \frac{1}{N}\big(S^{\frac{N}{2}}(\mu_1)+S^{\frac{N}{2}}(\mu_2)\big)$,  also a contradiction with (\ref{hhe6-5}). If $j_0=j_1=j$ and $j\not\in \{0, \infty\}$,   recall  (\ref{mla-e1}) and (\ref{mla-e2}) we see that
    $$c\geq \frac{1}{N}S\big(\rho_{j}^{\frac{2}{2^*}}+\overline{\rho}_{j}^{\frac{2}{2^*}}\big)\geq \frac{1}{N}\big(S^{\frac{N}{2}}(\mu_1)+S^{\frac{N}{2}}(\mu_2)\big),$$ also a contradiction with (\ref{hhe6-5}). Thus, we have
    either $u_n\rightarrow u_0$   or $v_n\rightarrow v_0$ strongly in $L^{2^*}(\R^N)$.

 \vskip 0.130in
 \noindent{\bf Step 3:} By the above arguments, under the conditions of $\displaystyle c<\frac{1}{N}\big(S^{\frac{N}{2}}(\mu_1)+S^{\frac{N}{2}}(\mu_2)\big)$ and (\ref{mla-e1}), we   obtain  that  either  $\rho_j$ or $\overline{\rho}_j$ equals  $0$ for any $j\not\in \{0,\infty\}$. Then (\ref{ml-e1}) is satisfied. Further,
   \be\lab{mla-e3}
   \rho_0\geq S^{\frac{N}{2}}(\mu_1), \rho_\infty\geq S^{\frac{N}{2}}(\mu_1), \overline{\rho}_0\geq S^{\frac{N}{2}}(\mu_2), \overline{\rho}_\infty\geq S^{\frac{N}{2}}(\mu_2)
   \ee
and $\mathcal{J}$  is finite.  For $j\in \mathcal{J}$ we have either
\begin{eqnarray}\lab{mla-e4}
 (\rho_j, \overline{\rho}_j)=(0, 0)\;\; \hbox{or}\;\big(\rho_j=0, \overline{\rho}_j\geq S^{\frac{N}{2}}\big)\;\hbox{or}\;\big(\rho_j\geq S^{\frac{N}{2}}, \overline{\rho}_j=0\big).
\end{eqnarray}
Thus, the steps 3 $\sim$  5 in Lemma \ref{hl-6-1} are valid here, and we finish the proof.
\ep

Consider $\displaystyle N=4,\alpha=\beta=2, \big(\frac{1-\mu_1}{1-\mu_2}\big)^{\frac{3}{2}}\geq \frac{1}{2}$. Similar to Lemma \ref{lemma-ml1}, we have the following lemma.
\bl\lab{lemma-amy1}
Consider $N=4, \alpha=\beta=2$. Assume $(H'_1),0\leq \mu_2\leq \mu_1<1,\mu_1+\mu_2\neq0, \big(\frac{1-\mu_1}{1-\mu_2}\big)^{\frac{3}{2}}\geq \frac{1}{2}$ and
\begin{eqnarray}\lab{lp-e1}
\Theta&\leq& \min\Big\{\frac{1-(1-\varepsilon_1)^{\frac{2}{N}}}{2^{\frac{\beta}{2}}\alpha (1-\varepsilon_1)^{\frac{\alpha-2}{2^*}}}, \frac{1-(1-\varepsilon_1)^{\frac{2}{N}}}{2^{\frac{\alpha}{2}}\beta (1-\varepsilon_1)^{\frac{\beta-2}{2^*}}},\frac{2-\sqrt{2}}{4}\Big\}\nonumber\\
&=&\min\{\frac{1-(1-\varepsilon_1)^{\frac{1}{2}}}{4}, \frac{2-\sqrt{2}}{4}\}\;\hbox{with $\varepsilon_1$ satisfying (\ref{ml-e13})}\nonumber\\
&=&\min\{\frac{2-(1-\mu_1)^{\frac{3}{2}}-(1-\mu_2)^{\frac{3}{2}}}{8}, \frac{2-\sqrt{2}}{4}\}.
\end{eqnarray}
Let $\{(u_n,v_n)\}\subset\overline{\mathcal{N}}$ be a Palais-Smale sequence for $\overline{\Phi}|_{\overline{\mathcal{N}}}$ at level $c\in\R$. Then, there exists a  constant $C$, such that  $||(u_n,v_n)||_{\D}\leq C$
for all $n\in \NN$ and that  $\overline{\Phi}'(u_n, v_n)\rightarrow 0$ in the dual space ${\D}^{*}$.  Moreover, if $c$ satisfies
\be\lab{amy-e-20}
\frac{1}{N}S^{\frac{N}{2}}(\mu_2)<c<\frac{1}{N}S^{\frac{N}{2}}(\mu_2)+\underset{(u, v)\in \mathcal{N}}{\inf}\Phi(u, v)\leq \frac{1}{N}\big(S^{\frac{N}{2}}(\mu_1)+S^{\frac{N}{2}}(\mu_2)\big);
\ee
\be\lab{amy-e-21}
c\neq \frac{l}{N}S^{\frac{N}{2}}(\mu_1)\;\; \hbox{and}\;\; c\neq \frac{l}{N}S^{\frac{N}{2}}\;\hbox{for all }l\in \NN\backslash\{0\},
\ee
then,  up to a subsequence, $(u_n,v_n)\rightarrow (u_0, v_0)$ in $\D$.
\el
\bp Follow  the processes of Lemma \ref{hl-6-1} and Lemma \ref{lemma-ml1} carefully. First,  by the similar argument as that in Lemma \ref{lemma-ml1} we obtain (\ref{mla-e3}) and (\ref{mla-e4}). Then the remaining work is similar to the proof of steps 3 $\sim$  5 in Lemma \ref{hl-6-1}. We can see that the only difference  is  that  from (\ref{amy-e-7}) to the end of Step 3. Thus, we can start as the following.  Since $N=4,\alpha=\beta=2$, we have
 $S(\mu)=(1-\mu)^{\frac{3}{4}}S,$ where $S$ is the best constant for the Sobolev inequality in $\R^4$.
Then (\ref{amy-e-4}) and (\ref{amy-e-7}) are equivalent to
\be\lab{amy-e-22}
\big(1-\mu_2\big)^{\frac{3}{4}}S\leq \|v_0\|_{L^4(\R^4)}^{2}+2 \Theta\big(1-\mu_2\big)^{\frac{3}{4}}S
\ee
and
\be\lab{amy-e-23}
\big(1-\mu_1\big)^{\frac{3}{4}}S\leq \|u_0\|_{L^4(\R^4)}^{2}+2 \Theta\big(1-\mu_2\big)^{\frac{3}{4}}S.
\ee
Consider $f(t)=t^{\frac{1}{2}}+2\Theta(1-\mu_2)^{\frac{3}{4}}S$, which is increasing  in $(0,+\infty)$. Assume   $\Theta\leq \frac{2-\sqrt{2}}{4},\mu_2\leq \mu_1$, we have  $ f\big(\frac{1}{2}S^{\frac{N}{2}}(\mu_2)\big) \leq  (1-\mu_1)^{\frac{3}{4}}S.$
Hence,  by (\ref{amy-e-23}), we have
$$\|u_0\|_{L^{2^*}(\R^N)}^{2^*}\geq \frac{1}{2}\big(1-\mu_2\big)^{\frac{3}{2}}S^2=\frac{1}{2}S^{\frac{N}{2}}(\mu_2).$$
Similarly, by (\ref{amy-e-22})  we have
$ \|v_0\|_{L^{2^*}(\R^N)}^{2^*}\geq \frac{1}{2}S^{\frac{N}{2}}(\mu_2).$
Thus
$$\|u_0\|_{L^{2^*}(\R^N)}^{2^*}+\|v_0\|_{L^{2^*}(\R^N)}^{2^*}\geq\frac{1}{2}S^{\frac{N}{2}}(\mu_2)
+\frac{1}{2}S^{\frac{N}{2}}(\mu_2)=S^{\frac{N}{2}}(\mu_2), $$
a contradiction to (\ref{61e-6}).
\ep

\vskip0.33in
\s{Nonexistence of the Nontrivial Least Energy Solution}

\renewcommand{\theequation}{4.\arabic{equation}}


We introduce the following  notation:

\begin{itemize}
 \item   [ ]  $C'_{\alpha,\beta}:=\frac{1}{\beta}\big(\frac{S(\mu_2)}{S(\mu_1)}-1\big)S^{\frac{4-(N-2)(\alpha+\beta-2)}{4}}(\mu_1)\;\;\hbox{if}\;\;\mu_2<\mu_1, \beta\geq 2;$
  \item  [ ] $ C'_{\alpha,\beta}:=\min\big\{\frac{2^{\frac{\beta-2}{2^*}}}{\beta},\frac{2^{\frac{\alpha-2}{2^*}}}{\alpha}\big\}
[1-(\frac{1}{2})^{\frac{2}{N}}]S^{\frac{4-(N-2)(\alpha+\beta-2)}{4}}(\mu)\;\;\hbox{if}\;\;  \mu_2=\mu_1=\mu $
\item   [  ]     \quad\quad\quad and $\alpha\geq2,\beta\geq2.$
   \end{itemize}
In particular, if $N=3, 4$, we define the following simpler constants：
  \begin{itemize}
 \item  [ ]  $C'_{\alpha,\beta}:=\min\{0.09, 0.09S^{-\frac{1}{2}}(\mu_1)\} $   if $N=3, \mu_2\leq \mu_1, \alpha\geq 2, \beta\geq 2$.
 \item  [ ]  $ C'_{\alpha, \beta}:=\frac{2-\sqrt{2}}{2}$ if  $N=4, \alpha=\beta=2, \mu_2\leq \mu_1$.

 \end{itemize}

\bt\lab{h-th5-1}
Assume that $  \beta\geq 2, \mu_2<\mu_1$ or $\alpha\geq 2, \beta\geq 2, \mu_2=\mu_1=\mu$.   Suppose that  $$ \begin{cases}&  \hbox{either }   h(x)\;\hbox{satisfies}\; (H_1)\;\hbox{if}\;\alpha+\beta<2^*;\\   & \hbox{or } h(x)\;\hbox{satisfies}\;(H'_1) \hbox{and }  \mu_1+\mu_2\neq 0\;\hbox{if}\;\alpha+\beta=2^*.  \end{cases}$$
If further, $\Theta\leq \Theta_0: =  C'_{\alpha,\beta}$  for respective  cases of  $\mu_1, \mu_2$ and $\alpha, \beta$, then the least energy  $c$ of the  system
$$c:=\underset{(u,v)\in\mathcal{N}}{\inf}\Phi(u, v)=\frac{1}{N}S^{\frac{N}{2}}(\mu_1).$$
Moreover, $c$ is achieved by and only by $(\pm z_{\sigma}^{\mu_1}, 0), \sigma>0$ if $\mu_2<\mu_1$ and by $(\pm z_\sigma^\mu, 0)$ and $(0, \pm z_\sigma^\mu)$ if $\mu_2=\mu_1=\mu\neq0$ \big(resp. $(\pm z_{\sigma,x_i}, 0)$ and $(0, \pm z_{\sigma,x_i})$ if $\mu_2=\mu_1=0$\big). That is,  problem (\ref{hhe6-4}) has no nontrivial least energy solution.
\et
\bp
First, we consider the case of $\mu_2<\mu_1$.
Since $(0, z_{\sigma}^{\mu_2})\in \mathcal{N}$ and $\Phi(0, z_{\sigma}^{\mu_2})=\frac{1}{N}S^{\frac{N}{2}}(\mu_2)$, we have that
$$c:=\underset{(u,v)\in \mathcal{N}}{\inf}\Phi(u, v)\leq \frac{1}{N}S^{\frac{N}{2}}(\mu_2).$$
On the other hand, note that $(z_{\sigma}^{\mu_1},0)\in \mathcal{N}$ and $\Phi(z_{\sigma}^{\mu_1}, 0)=\frac{1}{N}S^{\frac{N}{2}}(\mu_1)$, we get that
\be\lab{zou+c}c:=\underset{(u,v)\in \mathcal{N}}{\inf}\Phi(u, v)\leq \frac{1}{N}S^{\frac{N}{2}}(\mu_1).\ee
Since $\mu_2< \mu_1$, it follows that $S(\mu_1)< S(\mu_2)$ and that $c\leq \frac{1}{N}S^{\frac{N}{2}}(\mu_1).$
If $c<\frac{1}{N}S^{\frac{N}{2}}(\mu_1)$, then by Lemma \ref{group-lm1} we have $c>0$. By Lemma \ref{lm-6-3},   $c$ can be obtained by some $(0, 0)\not\equiv(\phi, \chi)\in \mathcal{N}$. Notice that the functional $\Phi$ is even, we have $(|\phi|, |\chi|)\in \mathcal{N}$  and $\Phi(\phi, \chi)=\Phi(|\phi|, |\chi|)$, hence $(|\phi|, |\chi|)$ is also a ground state for $\Phi$. Without loss of generality, we can assume $\phi\geq 0, \chi\geq 0$. Moreover, $\phi\not\equiv 0$ and $\chi\not\equiv 0$. If not, $\phi\equiv 0$ implies that $\chi\not\equiv 0$ is a solution of
$$\begin{cases} -\Delta v-\mu_2\frac{v}{|x|^2}=|v|^{2^*-2}v\quad \hbox{in}\; \R^N,\\
0\not\equiv v\in D^{1,2}(\R^N),\end{cases}$$
then $\Phi(\phi, \chi)=\Phi(0, \chi)\geq \frac{1}{N}S^{\frac{N}{2}}(\mu_2)> \frac{1}{N}S^{\frac{N}{2}}(\mu_1)>c$, a contradiction.
If $\chi\equiv 0, \phi\not\equiv 0$, then similarly we can get a contradiction.
Recalling that
$$c=\Phi(\phi, \chi)=\big(\frac{1}{2}-\frac{1}{\alpha+\beta}\big)\|(\phi, \chi)\|_{\D}^{2}+\big(\frac{1}{\alpha+\beta}-\frac{1}{2^*}\big)\Big(\|\phi\|_{L^{2^*}(\R^N)}^{2^*}+\|\chi\|_{L^{2^*}(\R^N)}^{2^*}\Big), $$
$$\int_{\R^N}\big(|\nabla \phi|^2-\mu_1\frac{\phi^2}{|x|^2}\big)dx =\int_{\R^N}|\phi|^{2^*}dx+\alpha\int_{\R^N} h(x)\phi^\alpha \chi^\beta dx,$$
$$\int_{\R^N} \big(|\nabla \chi|^2-\mu_2\frac{\chi^2}{|x|^2}\big)dx =\int_{\R^N}|\chi|^{2^*}dx+\beta\int_{\R^N} h(x)\phi^\alpha \chi^\beta dx. $$
Next, we  discuss two cases according to the sign of $\int_{\R^N} h(x)\phi^\alpha \chi^\beta dx.$
\begin{itemize}
\item Case 1: $\displaystyle\int_{\R^N} h(x)\phi^\alpha \chi^\beta dx\leq 0$.
\end{itemize}
It follows that  $$\int_{\R^N} \big(|\nabla \chi|^2-\mu_2\frac{\chi^2}{|x|^2}\big)dx\leq \int_{\R^N}|\chi|^{2^*}dx,$$
which implies that $ \|\chi\|_{L^{2^*}(\R^N)}^{2^*}\geq S^{\frac{N}{2}}(\mu_2).$
By the Hardy's inequality  and  the  Sobolev's  inequality, we have  that
\begin{align*}
\Phi(\phi, \chi)&=\big(\frac{1}{2}-\frac{1}{\alpha+\beta}\big)\|(\phi, \chi)\|_{\D}^{2}+\big(\frac{1}{\alpha+\beta}
-\frac{1}{2^*}\big)\Big(\|\phi\|_{L^{2^*}(\R^N)}^{2^*}+\|\chi\|_{L^{2^*}(\R^N)}^{2^*}\Big)\\
&\geq \big(\frac{1}{2}-\frac{1}{\alpha+\beta}\big)\|\chi\|_{\mu_2}^{2}+\big(\frac{1}{\alpha+\beta}
-\frac{1}{2^*}\big)\|\chi\|_{L^{2^*}(\R^N)}^{2^*}\\
&\geq \big(\frac{1}{2}-\frac{1}{\alpha+\beta}\big)S(\mu_2)\|\chi\|_{L^{2^*}(\R^N)}^{2}+\big(\frac{1}{\alpha+\beta}
-\frac{1}{2^*}\big)\|\chi\|_{L^{2^*}(\R^N)}^{2^*}\\
&\geq \frac{1}{N}S^{\frac{N}{2}}(\mu_2)>\frac{1}{N}S^{\frac{N}{2}}(\mu_1)>c,
\end{align*}
a contradiction.
\begin{itemize}
\item Case 2: $\displaystyle\int_{\R^N} h(x)\phi^\alpha \chi^\beta dx> 0$.
\end{itemize}
Notice that
\begin{eqnarray}\lab{he-5-1}
c=\Phi(\phi, \chi)&=&\frac{1}{N}\big(\|\phi\|_{L^{2^*}(\R^N)}^{2^*}+\|\chi\|_{L^{2^*}(\R^N)}^{2^*}\big)\nonumber\\
&&+(\alpha+\beta)(\frac{1}{2}-\frac{1}{\alpha+\beta})\int_{\R^N} h(x)\phi^\alpha \chi^\beta dx\nonumber\\
&<&\frac{1}{N}S^{\frac{N}{2}}(\mu_1) ,
\end{eqnarray}
we have
\be\lab{hsn-e3}\|\chi\|_{L^{2^*}(\R^N)}^{2^*}<S^{\frac{N}{2}}(\mu_1), \|\phi\|_{L^{2^*}(\R^N)}^{2}<S^{\frac{N}{2}}(\mu_1).\ee
By  the Hardy's  inequality (or the Sobolev's inequality) and  the H\"{o}lder's inequality, we see that
\begin{eqnarray}\lab{amy-e-3}
S(\mu_2)\|\chi\|_{L^{2^*}(\R^N)}^{2}&\leq& \|\chi\|_{\mu_2}^{2}\nonumber\\
&=&\|\chi\|_{L^{2^*}(\R^N)}^{2^*}+\beta \int_{\R^N} h(x)\phi^\alpha \chi^\beta dx\nonumber\\
&\leq &\|\chi\|_{L^{2^*}(\R^N)}^{2^*}+ \beta \int_{\R^N} h_+(x)\phi^\alpha \chi^\beta dx\nonumber\\
&\leq &\|\chi\|_{L^{2^*}(\R^N)}^{2^*}+ \beta \Theta S^{\frac{N-2}{4}\alpha}(\mu_1)\|\chi\|_{L^{2^*}(\R^N)}^{\beta}.
\end{eqnarray}
Denote $\displaystyle\sigma=\|\chi\|_{L^{2^*}(\R^N)}^{2^*}$, then
\be\lab{hsn-e2}S(\mu_2)\leq \sigma^{1-\frac{2}{2^*}}+\beta\Theta S^{\frac{N-2}{4}\alpha}(\mu_1)\sigma^{\frac{\beta-2}{2^*}}.\ee
Let $f(\sigma)=\sigma^{1-\frac{2}{2^*}}+\beta\Theta S^{\frac{N-2}{4}\alpha}(\mu_1)\sigma^{\frac{\beta-2}{2^*}}$.  Since $\beta\geq 2$, we have that $f(\sigma)$ is increasing in $(0, +\infty)$.
If $\mu_2<\mu_1, \beta\geq 2, \Theta\leq C'_{\alpha,\beta}$, the  direct calculation implies that  $
f\big(S^{\frac{N}{2}}(\mu_1)\big) \leq  S(\mu_2).$
Then (\ref{hsn-e2}) implies that $\displaystyle \|\chi\|_{L^{2^*}(\R^N)}^{2^*}\geq S^{\frac{N}{2}}(\mu_1)$, a contradiction to (\ref{hsn-e3}).
By the above arguments, we obtain that $\displaystyle c=\frac{1}{N}S^{\frac{N}{2}}(\mu_1)$ and obviously it is achieved by $(\pm z_{\sigma}^{\mu_1}, 0), \sigma>0$ .

Let $(\phi, \chi)$ be a minimizer  of  $c=\frac{1}{N}S^{\frac{N}{2}}(\mu_1)$.  If $\phi\not\equiv 0$ and $\chi\not\equiv 0$, it will lead to a contradiction if repeating the above arguments.
Here, we may assume $\phi\equiv 0$, then $\Phi(\phi, \chi)\geq \frac{1}{N}S^{\frac{N}{2}}(\mu_2)> \frac{1}{N}S^{\frac{N}{2}}(\mu_1)=c$ if $\mu_2<\mu_1$. Hence, $c$ is only achieved by $(\phi, 0)$, where $\phi$ is a  weak solution of
$$\begin{cases}-\Delta u-\mu_1\frac{u}{|x|^2}=|u|^{2^*-2}u\quad \hbox{in}\;\R^N,\\
0\not\equiv u\in D^{1,2}(\R^N).\end{cases}$$
By \cite{T.S}, it is easy to see that $c$ is only achieved by $(\pm z_\sigma^\mu, 0), \sigma>0$.

 Second, we consider the case of $\mu_2=\mu_1=\mu, \alpha\geq 2, \beta\geq 2, \Theta\leq C'_{\alpha, \beta}$.  The difference is the computation of $f\big(\frac{1}{2}S^{\frac{N}{2}}(\mu_1)\big)$:
\begin{eqnarray*}
f\big(\frac{1}{2}S^{\frac{N}{2}}(\mu_1)\big)&=&\big[\frac{1}{2}S^{\frac{N}{2}}(\mu_1)\big]^{\frac{2}{N}}
+\beta\Theta S^{\frac{N-2}{4}\alpha}(\mu_1)\big[\frac{1}{2}S^{\frac{N}{2}}(\mu_1)\big]^{\frac{\beta-2}{2^*}}\\
&=&(\frac{1}{2})^{\frac{2}{N}}S(\mu_1)+\beta\Theta(\frac{1}{2})^{\frac{\beta-2}{2^*}}S^{\frac{(N-2)(\alpha+\beta-2)}{4}}(\mu_1)\\
&\leq &S(\mu_2)=S(\mu_1).
\end{eqnarray*}
Then (\ref{hsn-e2}) implies that $\displaystyle \|\chi\|_{L^{2^*}(\R^N)}^{2^*}\geq \frac{1}{2}S^{\frac{N}{2}}(\mu)$. Similarly, we can obtain that
$$\|\phi\|_{L^{2^*}(\R^N)}^{2^*}\geq \frac{1}{2}S^{\frac{N}{2}}(\mu).$$
Then
$$\|\chi\|_{L^{2^*}(\R^N)}^{2^*}+\|\phi\|_{L^{2^*}(\R^N)}^{2^*}\geq S^{\frac{N}{2}}(\mu)=S^{\frac{N}{2}}(\mu_1),$$
a contradiction to (\ref{he-5-1}). Finally, it is easy to  show  that $c$ is achieved by and only by semi-trivial solutions $(\pm z_\sigma^\mu, 0)$ and $(0, \pm z_\sigma^\mu)$ if  $\mu\neq 0$ (resp.$(\pm z_{\sigma,x_i}, 0)$ and $(0, \pm z_{\sigma,x_i})$  if  $\mu= 0$  ).
\ep

\noindent{\bf  Proof of  Theorem \ref{h-th5-1-1}}. It is a straightforward consequence of Theorem  \ref{h-th5-1}.  $\Box$


\vskip0.33in
\s{The Existence of Mountain Pass Solutions}

\renewcommand{\theequation}{5.\arabic{equation}}

If $\min\{\alpha, \beta\}=2$, we  denote
\be\label{zwm=may-2}  C_{\alpha,\beta,\mu_1,\mu_2}:=\min\Big\{\frac{S(\mu_1)}{2S^{\frac{N-2}{4}\beta}(\mu_2)},
    \frac{S(\mu_2)}{2S^{\frac{N-2}{4}\alpha}(\mu_1)} \Big\}\ee


\br\lab{zou+1} We observe that
     $C_{\alpha,\beta,\mu_1, \mu_2}>\min\{0.3, 0.3S^{\frac{1}{2}}(\mu_1)\}   \hbox{ if}\; N=3;$
     $C_{\alpha,\beta,\mu_1, \mu_2}>\frac{\sqrt{2}}{4}\;\hbox{if}\; N=4$ and $ \alpha=\beta=2.$

\er

\br\lab{hr6-1}
For $\alpha>2$, we can see that the indefinite sign of $h(x)$  has no effect on the proof of Theorem 2.2 in  \cite{BVI}. Thus,  for all $\sigma>0, z_{\sigma}^{\mu_2}$ is a local minimum point of $\Phi$ in $\mathcal{N}$.
Denote
 $$\mathcal{N}_{\mu_2}=\big\{u\in D^{1,2}(\R^N)\backslash \{0\}:\|u\|_{\mu_2}^{2}=\|u\|_{L^{2^*}(\R^N)}^{2^*}\big\},$$ which  is  the corresponding Nehari manifold of $I_{\mu_2}$,
 where $I_\mu$ is defined in (\ref{hhe6-2}).
For $\alpha=2, \mu_2>0$, we can obtain that \big(for the details we refer to \cite[Theorem 2.2]{BVI} or \cite[Lemma 3.4]{ZZ}\big)  $(\phi, z_{\sigma}^{\mu_2}+\chi)\in\mathcal{N}$, where $z_\sigma^\mu$ is defined in (\ref{zzz-1}).  Let $t>0$ satisfy  $t(z_{\sigma}^{\mu_2}+\chi)\in \mathcal{N}_{\mu_2}$, then
\begin{eqnarray*}
&&\Phi(\phi, z_{\sigma}^{\mu_2}+\chi)-\Phi(0, t(z_{\sigma}^{\mu_2}+\chi))\\
& &=\big(\frac{1}{2}\|\phi\|_{\mu_1}^{2}-\int_{\R^N}h(x)|\phi|^2|z_{\sigma}^{\mu_2}+\chi|^\beta\big)(1+o(1))\\
&&\geq\Big(\frac{1}{2}S(\mu_1)\|\phi\|_{L^{2^*}(\R^N)}^{2}-\Theta\|\phi\|_{L^{2^*}(\R^N)}^{2}\|z_{\sigma}^{\mu_2}+\chi\|_{L^{2^*}(\R^N)}^{\beta}\Big)(1+o(1))\\
& &=\Big(\frac{1}{2}S(\mu_1)-\Theta\|z_{\sigma}^{\mu_2}+\chi\|_{L^{2^*}(\R^N)}^{\beta}\Big)\|\phi\|_{L^{2^*}(\R^N)}^{2}(1+o(1))\\
& &\rightarrow\Big(\frac{1}{2}S(\mu_1)-\Theta S^{\frac{N-2}{4}\beta}(\mu_2)\Big)\|\phi\|_{L^{2^*}(\R^N)}^{2}(1+o(1))
\end{eqnarray*}
 as $\|(\phi,\chi)\|_{\D}\rightarrow 0$. Thus, when   $\Theta<  C_{\alpha,\beta,\mu_1,\mu_2}$ (see  (\ref{zwm=may-2} ) ), we have that
 $$\frac{1}{2}S(\mu_1)-\Theta S^{\frac{N-2}{4}\beta}(\mu_2)>0.$$
It follows that
\be\lab{ze6-1}
\Phi(\phi, z_{\sigma}^{\mu_2}+\chi)-\Phi(0, t(z_{\sigma}^{\mu_2}+\chi))>0
\ee
 provided that  $(\phi, z_{\sigma}^{\mu_2}+\chi)\in \mathcal{N}$  is sufficiently closed to $(0, z_{\sigma}^{\mu_2})$.
  On the other hand, since $z_{\sigma}^{\mu_2}$ is a minimizer of $I_{\mu_2}$ on $\mathcal{N}_{\mu_2}$, we have
  \be\lab{ze6-2}
  \Phi(0, t(z_{\sigma}^{\mu_2}+\chi)-\Phi(0, z_{\sigma}^{\mu_2}))=I_{\mu_2}\big(t(z_{\sigma}^{\mu_2}+\chi)\big)-I_{\mu_2}(z_{\sigma}^{\mu_2})\geq 0.
  \ee
 From (\ref{ze6-1}) and (\ref{ze6-2}), we conclude that
  $\Phi(\phi, z_{\sigma}^{\mu_2}+\chi)-\Phi(0, z_{\sigma}^{\mu_2})\geq 0$
 for any $(\phi, z_{\sigma}^{\mu_2}+\chi)\in \mathcal{N}$ sufficiently close  to $(0, z_{\sigma}^{\mu_2})$(with strictly inequality hold outside the manifold $\{0\}\times Z_{\mu_2}$, where $Z_{\mu_2}=\{z_{\sigma}^{\mu_2}, \sigma>0\}$ ),
i.e., $(0, z_{\sigma}^{\mu_2})$ is a local minimum point of $\Phi$ in $\mathcal{N}$.
Similarly, for $\alpha=2, \mu_2=0$,  we replace $z_{\sigma}^{\mu_2}$ by $z_{\sigma, x_i}$, where $z_{\sigma, x_i}=\sigma^{-\frac{N-2}{2}}z_1(\frac{x-x_i}{\sigma}), \sigma>0, x_i\in\R^N$, and $z_\sigma(x)=z_{\sigma, 0}$, $\displaystyle z_{1,0}(x)=\frac{[N(N-2)]^{\frac{N-2}{4}}}{[1+|x|^2]^{\frac{N-2}{2}}}$, we can obtain that
 $\Phi(\phi, z_{\sigma,x_i}+\chi)-\Phi(0, z_{\sigma, x_i})\geq 0$
for any $(\phi, z_{\sigma, x_i}+\chi)\in \mathcal{N}$ sufficiently closed to $(0, Z_{\sigma, x_i})$ (with strictly inequality holding outside the manifold $\{0\}\times Z_0$, where $Z_0=\{z_{\sigma, x_i}, \sigma>0, x_i\in \R^N\}$), i.e., $(0, z_{\sigma, x_i})$ is a local minimum point of $\Phi$ in $\mathcal{N}$.
Similarly, if $\beta>2$ or $\beta=2$ with  $\Theta<  C_{\alpha,\beta,\mu_1,\mu_2}$ (see  (\ref{zwm=may-2} ) ), we can obtain that $(z_{\sigma}^{\mu_1}, 0)$ is a local minimum. So does $(z_{\sigma,x_i},0)$ if $\mu_1=0$ and  either $ \beta>2$ or $\beta=2 $ with    $\Theta<  C_{\alpha,\beta,\mu_1,\mu_2}$. \qed
\er


\vskip0.23in
By Remark \ref{hr6-1},   the semi-trivial solutions $(0, z_{\sigma}^{\mu_2})$ or $(0, z_{\sigma, x_i})$ turns  out to be local minimum points for the functional $\Phi\big|_{\mathcal{N}}$, which consequently exhibits a mountain pass geometry.

\vskip0.1336in
Next, our goal is to construct a mountain pass level for the functional on the Nehari manifold at which the Palais-Smale condition holds in view of  Lemmas \ref{hl-6-1},   \ref{lemma-ml1} and  \ref{lemma-amy1}. For  the simplicity, when $\mu=0$, we also use the  notation $z_1^0$ instead of $z_{1,0}$. When $\alpha>2$ or $\alpha=2$ with  $\Theta<  C_{\alpha,\beta,\mu_1,\mu_2}$ (see  (\ref{zwm=may-2} ) ),  by Remark \ref{hr6-1}, we know  that $(z_{1}^{\mu_1},0)$ is a local minimum. Similarly, when $\beta>2$ or $\beta=2$ with
 $\Theta<  C_{\alpha,\beta,\mu_1,\mu_2}$, $(0, z_{1}^{\mu_2})$ is a local minimum. By Theorem \ref{h-th5-1}, when $\displaystyle \mu_2<\mu_1, \Theta\leq C'_{\alpha,\beta}$ or $\displaystyle\mu_2=\mu_1, \Theta\leq C'_{\alpha,\beta}, $ we have $ \underset{(u,v)\in\mathcal{N}}{\inf}\Phi(u,v)=\frac{1}{N}S^{\frac{N}{2}}(\mu_1)$. Since $\Phi$ is even, for $(u, v)\in \mathcal{N}$, we have $(|u|, |v|)\in\mathcal{N}$. Thus, $(|u|, |v|)\in \mathcal{N}\cap\overline{\mathcal{N}}$ and $\overline{\Phi}(|u|, |v|)=\Phi(|u|, |v|)=\Phi(u, v)$, then $\displaystyle \underset{(u, v)\in \overline{\mathcal{N}}}{\inf}\overline{\Phi}(u,v)\leq \frac{1}{N}S^{\frac{N}{2}}(\mu_1)$.  Assume  $\displaystyle \underset{(u, v)\in \overline{\mathcal{N}}}{\inf}\overline{\Phi}(u,v)< \frac{1}{N}S^{\frac{N}{2}}(\mu_1)\leq \frac{1}{N}S^{\frac{N}{2}}(\mu_2)$, similar to Lemma \ref{lm-6-3},  we can prove that the infimum  can be achieved by some $(u_0, v_0)$, and the minimizer is a critical point of $\overline{\Phi}$. It is easy to see $u_0\geq 0, v_0\geq 0$, thus $(u_0, v_0)\in \mathcal{N}\cap \overline{\mathcal{N}}$ and
$\overline{\Phi}(u_0, v_0)=\Phi(u_0, v_0)\geq \frac{1}{N}S^{\frac{N}{2}}(\mu_1)$, a contradiction. Thus, we have
\be\lab{bb-e5}
 \underset{(u,v)\in\overline{\mathcal{N}}}{\inf}\overline{\Phi}(u,v)=\frac{1}{N}S^{\frac{N}{2}}(\mu_1)
\ee
 when  $ \mu_2\leq \mu_1,    \Theta\leq C'_{\alpha,\beta}.$

 \vskip0.1in

Next,   we consider the set of the paths in $\overline{\mathcal{N}}$ joining $(z_{\sigma}^{\mu_1}, 0)$ with $(0, z_{\sigma}^{\mu_2})$, i.e.,
$$\Sigma=\left\{\begin{array}{lll}
&\gamma=(\gamma_1, \gamma_2):[0,1]\rightarrow \overline{\mathcal{N}}\;\hbox{continous such that}\\\\
&\gamma_1(0)=0,\gamma_1(1)=z_{\sigma}^{\mu_1},\gamma_2(0)=z_{\sigma}^{\mu_2},\gamma_2(1)=0\end{array}\right\}
,$$
and define the associated mountain pass level
\be\lab{nn-e2}
C_{MP}:=\underset{\gamma\in\Sigma}{\inf}\underset{t\in[0,1]}{\max}\overline{\Phi}(\gamma(t)).
\ee

\vskip0.336in

\subsection{The case of $N=3, \frac{1}{2}<\frac{1-4\mu_1}{1-4\mu_2}, S^{\frac{N}{2}}(\mu_2)+S^{\frac{N}{2}}(\mu_1)\leq S^{\frac{N}{2}}$.}
If $\mu_2<\mu_1$, we define
\be\lab{xiu-e1}
M_1:=1-\frac{2\sqrt{2}(1-4\mu_1)}{[(1-4\mu_2)^{\frac{2}{3}}+(1-4\mu_1)^{\frac{2}{3}}]^{\frac{3}{2}}},
\ee
\be\lab{xiu-e2}
M_2:=\frac{2(1-4\mu_1)}{1-4\mu_2}-1,
\ee
\be\lab{xiu-e3}
M_3:=\frac{1}{2}\min\{M_1,M_2\},
\ee
then we have
\be\lab{xiu-e4}
\frac{2}{N}(1-M_3)\big(\frac{S(\mu_2)+S(\mu_1)}{2}\big)^{\frac{N}{2}}
>\frac{2}{N}S^{\frac{N}{2}}(\mu_1)>\frac{1+M_3}{N}S^{\frac{N}{2}}(\mu_2).
\ee
If $\mu_2=\mu_1=\mu$, we define
\be\lab{lp-ez1}
M_3:=\frac{1}{4},
\ee
then we have
\be\lab{lp-ez2}
\frac{2}{N}(1-M_3)\big(\frac{S(\mu_1)+S(\mu_2)}{2}\big)^{\frac{N}{2}}>\frac{1+M_3}{N}S^{\frac{N}{2}}(\mu_2).
\ee
Define
\be\lab{xiu-e5}
M_4:=\frac{1-\frac{1}{2}(1-M_3)^{\frac{2}{N}}}{(\alpha+\beta)(\frac{1}{2})^{\frac{(N-2)(\alpha+\beta-2)}{4}}}>\frac{\sqrt{2}}{12}.
\ee

\bl\lab{lp-l1}
Assume $N=3, \alpha\geq 2, \beta\geq 2, \alpha+\beta\leq 2^*,   \frac{1}{2}<\frac{1-4\mu_1}{1-4\mu_2}$  such that  $\displaystyle S^{\frac{N}{2}}(\mu_1)+S^{\frac{N}{2}}(\mu_2)\leq S^{\frac{N}{2}}$.   The weight function $h(x)$ satisfies $(H_2)$ and $\begin{cases} (H_1)\;\hbox{if}\;\alpha+\beta<2^*\\ (H'_1)\;\hbox{if}\;\alpha+\beta=2^* \end{cases}$. Moreover, if $\min\{\alpha, \beta\}=2$, we   assume  that   $\Theta<  C_{\alpha,\beta,\mu_1,\mu_2}$ (see  (\ref{zwm=may-2} ) ), where $\Theta$ is defined in (\ref{ha1}).  Then, if $\Theta\leq \min\{M_4, C'_{\alpha,\beta}\}$,  the functional $\overline{\Phi}$ exhibits a mountain pass geometry and the mountain pass level satisfies (\ref{hhe6-5}) and (\ref{hhe6-6}).
\el

\bp We claim that when $\widetilde{\alpha}\leq M_4$, we have
\be\lab{mm-e1}
\underset{t\in[0,1]}{\max}\overline{\Phi}\big(\gamma(t)\big)\geq \frac{2}{N}(1-M_3)\Big(\frac{S(\mu_1)+S(\mu_2)}{2}\Big)^{\frac{N}{2}} \;\hbox{for all}\;\gamma\in \Sigma.
\ee
Let $(\gamma_1, \gamma_2)\in \Sigma$, since $(\gamma_1(t), \gamma_2(t))\in \overline{\mathcal{N}}$, we have
\begin{eqnarray*}
&&\int_{\R^N}|\nabla \gamma_1(t)|^2dx-\mu_1\int_{\R^N}\frac{\gamma_1^2(t)}{|x|^2}dx+\int_{\R^N}|\nabla \gamma_2(t)|^2dx-\mu_2\int_{\R^N}\frac{\gamma_2^2(t)}{|x|^2}dx\\
& &=\int_{\R^N}\big(\gamma_1(t)\big)_{+}^{2^*}dx+\int_{\R^N}\big(\gamma_2(t)\big)_{+}^{2^*}dx+
(\alpha+\beta)\int_{\R^N}h(x)\big(\gamma_1(t)\big)_+^\alpha\big(\gamma_2(t)\big)_+^\beta dx
\end{eqnarray*}
and
\begin{eqnarray}\lab{lap-e1}
\overline{\Phi}\big(\gamma_1(t), \gamma_2(t)\big)&=&\frac{1}{N}\Big(\int_{\R^N}\big(\gamma_1(t)\big)_{+}^{2^*}dx+\int_{\R^N}\big(\gamma_2(t)\big)_{+}^{2^*}dx\Big)\nonumber\\
&&+\frac{\alpha+\beta-2}{2}\int_{\R^N}h(x)\big(\gamma_1(t)\big)_+^\alpha\big(\gamma_2(t)\big)_+^\beta dx
\end{eqnarray}
or
\begin{eqnarray}\lab{lap-e2}
\overline{\Phi}\big(\gamma_1(t), \gamma_2(t)\big)&=&\big(\frac{1}{2}-\frac{1}{\alpha+\beta}\big)\|(\gamma_1(t), \gamma_2(t))\|_{\D}^{2}+\big(\frac{1}{\alpha+\beta}-\frac{1}{2^*}\big)\nonumber\\
&&\Big(\int_{\R^N}\big(\gamma_1(t)\big)_{+}^{2^*}dx+\int_{\R^N}\big(\gamma_2(t)\big)_{+}^{2^*}dx\Big)
\end{eqnarray}
or
\begin{eqnarray}\lab{lap-e3}
\overline{\Phi}\big(\gamma_1(t), \gamma_2(t)\big)&=&\frac{1}{N}\|(\gamma_1(t), \gamma_2(t))\|_{\D}^{2}-\big(1-\frac{\alpha+\beta}{2^*}\big)\nonumber\\
&&\int_{\R^N}h(x)\big(\gamma_1(t)\big)_+^\alpha\big(\gamma_2(t)\big)_+^\beta dx.
\end{eqnarray}
Denote $f_i(t)=\int_{\R^N}\big(\gamma_i(t)\big)_{+}^{2^*}dx$ for $i=1,2$, then
$$f_1(0)=0<f_2(0)=\int_{\R^N}(z_{\sigma}^{\mu_2})^{2^*}dx\;\hbox{and}\;f_1(1)=\int_{\R^N}(z_{\sigma}^{\mu_1})^{2^*}dx>f_2(1)=0,$$
hence, by  the  continuity, there exists some $t_0\in (0,1)$ such that $f_1(t_0)=f_2(t_0)>0$.
From the definition of $S(\mu_1)$ and $S(\mu_2)$ and the H\"{o}lder's  inequality, we have
\begin{eqnarray*}
&&S(\mu_1)\Big(\int_{\R^N}\big(\gamma_1(t_0)\big)_{+}^{2^*}\Big)^{\frac{2}{2^*}}+S(\mu_2)\Big(\int_{\R^N}\big(\gamma_2(t_0)\big)_{+}^{2^*}\Big)^{\frac{2}{2^*}}\\
& & \leq \int_{\R^N}|\nabla \gamma_1(t_0)|^2dx-\mu_1\int_{\R^N}\frac{\gamma_1^2(t_0)}{|x|^2}dx+\int_{\R^N}|\nabla \gamma_2(t_0)|^2dx-\\
& & \quad \mu_2\int_{\R^N}\frac{\gamma_2^2(t_0)}{|x|^2}dx\\
& &=\int_{\R^N}\big(\gamma_1(t_0)\big)_{+}^{2^*}dx+\int_{\R^N}\big(\gamma_2(t_0)\big)_{+}^{2^*}dx+\\
& & \quad
(\alpha+\beta)\int_{\R^N}h(x)\big(\gamma_1(t_0)\big)_+^\alpha\big(\gamma_2(t_0)\big)_+^\beta dx\\
& &\leq \int_{\R^N}\big(\gamma_1(t_0)\big)_{+}^{2^*}dx+\int_{\R^N}\big(\gamma_2(t_0)\big)_{+}^{2^*}dx+\\
&&\quad (\alpha+\beta)\Theta \Big(\int_{\R^N}\big(\gamma_1(t_0)\big)_{+}^{2^*}dx\Big)^{\frac{\alpha}{2^*}}
\Big(\int_{\R^N}\big(\gamma_2(t_0)\big)_{+}^{2^*}dx\Big)^{\frac{\beta}{2^*}}.
\end{eqnarray*}
 Set $\sigma=\int_{\R^N}\big(\gamma_1(t_0)\big)_{+}^{2^*}dx=\int_{\R^N}\big(\gamma_2(t_0)\big)_{+}^{2^*}dx$, then we obtain that
 \be\lab{mm-e2}
 \big(S(\mu_1)+S(\mu_2)\big)\sigma^{\frac{2}{2^*}}\leq \sigma+(\alpha+\beta)\Theta\sigma^{\frac{\alpha+\beta}{2^*}}
 \ee
and it is equivalent to
\be\lab{mm-e3}
S(\mu_1)+S(\mu_2)\leq \sigma^{1-\frac{2}{2^*}}+(\alpha+\beta)\Theta\sigma^{\frac{\alpha+\beta-2}{2^*}}.
\ee
Denote $g(\sigma):=\sigma^{1-\frac{2}{2^*}}+(\alpha+\beta)\Theta\sigma^{\frac{\alpha+\beta-2}{2^*}}$,  which is increasing  in $(0,\infty)$. Recall that $\Theta\leq M_4$, it follows that
\begin{eqnarray*}
&&g\Big((1-M_3)\big(\frac{S(\mu_1)+S(\mu_2)}{2}\big)^{\frac{N}{2}}\Big)\\
& & = (1-M_3)^{\frac{2}{N}}\frac{S(\mu_1)+S(\mu_2)}{2}+(\alpha+\beta)\Theta\Big(\frac{S(\mu_1)+S(\mu_2)}{2}\Big)^{\frac{(\alpha+\beta-2)(N-2)}{4}}\\
& &\leq  S(\mu_1)+S(\mu_2),
\end{eqnarray*}
which implies that
$$\sigma\geq (1-M_3)(\frac{S(\mu_1)+S(\mu_2)}{2})^{\frac{N}{2}},$$
then
\begin{eqnarray*}
\overline{\Phi}(\gamma(t_0))&\geq & 2\Big[(\frac{1}{2}-\frac{1}{\alpha+\beta})(\frac{S(\mu_1)+S(\mu_2)}{2})\sigma^{\frac{2}{2^*}}+(\frac{1}{\alpha+\beta}-\frac{1}{2^*})\sigma\Big]\\
&\geq&2\Big[(\frac{1}{2}-\frac{1}{\alpha+\beta})(1-M_3)^{\frac{2}{2^*}}(\frac{S(\mu_1)+S(\mu_2)}{2})^{\frac{N}{2}}\\
&&+(\frac{1}{\alpha+\beta}-\frac{1}{2^*})(1-M_3)(\frac{S(\mu_1)+S(\mu_2)}{2})^{\frac{N}{2}}\Big]\\
&>&\frac{2}{N}(1-M_3)(\frac{S(\mu_1)+S(\mu_2)}{2})^{\frac{N}{2}},
\end{eqnarray*}
we prove the claim (\ref{mm-e1}) and   obtain that
$$C_{MP}>\frac{1+M_3}{N}S^{\frac{N}{2}}(\mu_2)=(1+M_3)\overline{\Phi}(0, z_{\sigma}^{\mu_2}),$$ and hence $\overline{\Phi}$ exhibits a mountain pass geometry. In particular,  $C_{MP}>\frac{2}{N}S^{\frac{N}{2}}(\mu_1)$ if $\mu_2<\mu_1$;    $C_{MP}>\frac{1+M_3}{N}S^{\frac{N}{2}}(\mu_1)$ if $\mu_2=\mu_1$.

\vskip0.12in

Next, we consider a special path  $\gamma(t)=k(t)\big(\gamma_1(t), \gamma_2(t)\big)\in \Sigma, t\in[0,1]$ with $\gamma_1(t)=t^{\frac{1}{2}}z_{\sigma}^{\mu_1}$ and $\gamma_2(t)=(1-t)^{\frac{1}{2}}z_{\sigma}^{\mu_2}$, where $k(t)$ is a positive function such that $k(t)\big(\gamma_1(t), \gamma_2(t)\big)\in \mathcal{N}\cap \overline{\mathcal{N}}$. By the definition of the Nehari manifold, $k(t)$ is well defined and unique. For  the simplicity, we   set
$$a:=\|z_{\sigma}^{\mu_1}\|_{\mu_1}^{2}=\int_{\R^N}|z_{\sigma}^{\mu_1}|^{2^*}dx=S^{\frac{N}{2}}(\mu_1)$$
and
$$b:=\|z_{\sigma}^{\mu_2}\|_{\mu_2}^{2}
=\int_{\R^N}|z_{\sigma}^{\mu_2}|^{2^*}dx=S^{\frac{N}{2}}(\mu_2).$$
Since $k(t)\big(\gamma_1(t), \gamma_2(t)\big)\in \mathcal{N}\cap \overline{\mathcal{N}}$, then by $(H_2)$, we obtain
\begin{eqnarray*}
&& \|(t^{\frac{1}{2}}z_{\sigma}^{\mu_1}, (1-t)^{\frac{1}{2}}z_{\sigma}^{\mu_2})\|_{\D}^{2}\\
&&=k^{2^*-2}(t)\big((1-t)^{\frac{2^*}{2}}a+t^{\frac{2^*}{2}}b\big) \\
& & \quad +(\alpha+\beta)k^{\alpha+\beta-2}(t)(1-t)^{\frac{\alpha}{2}}t^{\frac{\beta}{2}}\int_{\R^N}h(x)|z_{\sigma}^{\mu_1}|^\alpha |z_{\sigma}^{\mu_2}|^\beta dx\\
&& >   k^{2^*-2}(t)\big((1-t)^{\frac{2^*}{2}}a+t^{\frac{2^*}{2}}b\big) \;\hbox{for}\; 0<t<1.
\end{eqnarray*}
Hence,
\be\lab{bb-e1}
k(t)<\Big[\frac{(1-t)a+tb}{(1-t)^{\frac{2^*}{2}}a+t^{\frac{2^*}{2}}b}\Big]^{\frac{N-2}{4}}\;\hbox{for all}\; 0<t<1
\ee
and $k(0)=k(1)=1$.
Combine  $(H_2)$ and (\ref{lap-e3}), it follows that
$$\overline{\Phi}\big(k(t)(\gamma_1(t), \gamma_2(t))\big)<\frac{k^2(t)}{N}\big((1-t)a+tb\big) \;\hbox{for all}\; 0<t<1,$$
hence,
\begin{eqnarray}\lab{bb-e2}
C_{MP} = \underset{\gamma\in \Sigma}{\inf}\underset{t\in[0,1]}{\max}\overline{\Phi}\big(\gamma(t)\big)
 \leq \underset{t\in[0,1]}{\max}\overline{\Phi}(k(t)(\gamma_1(t), \gamma_2(t))).
\end{eqnarray}
By (\ref{lap-e3}) and (\ref{bb-e1}), we have
$$\overline{\Phi}\big(k(t)(\gamma_1(t), \gamma_2(t))\big)<\frac{1}{N}\Big[\frac{(1-t)a+tb}{(1-t)^{\frac{2^*}{2}}a+t^{\frac{2^*}{2}}b}\Big]^{\frac{N-2}{N}}\big((1-t)a+tb\big)\;\hbox{for all}\; t\in (0,1).$$
After a direct computation, the right-hand side achieves its maximum at $t=\frac{1}{2}$ and the maximum is $\frac{1}{N}(a+b)$. Combining with (\ref{bb-e2}), we obtain that
\be\lab{bb-e3}
C_{MP}\leq \underset{t\in[0,1]}{\max}\overline{\Phi}(k(t)(\gamma_1(t), \gamma_2(t)))<\frac{1}{N}\big(S^{\frac{N}{2}}(\mu_1)+S^{\frac{N}{2}}(\mu_2)\big).
\ee
Recall that $\frac{1}{2}<\frac{1-4\mu_1}{1-4\mu_2}$, we have $S^{\frac{N}{2}}(\mu_2)<2 S^{\frac{N}{2}}(\mu_1)$. Therefore, from  the above arguments we obtain that
$$\frac{1+M_3}{2}S^{\frac{N}{2}}(\mu_2)<\frac{2}{N}S^{\frac{N}{2}}(\mu_1)<C_{MP}<\frac{1}{N}(S^{\frac{N}{2}}(\mu_1)+S^{\frac{N}{2}}(\mu_2))
$$
$$<\min\{\frac{1}{N}S^{\frac{N}{2}}, \frac{3}{N}S^{\frac{N}{2}}(\mu_1)\}  \;\hbox{if}\;\mu_2<\mu_1; $$
and
$$\frac{1+M_3}{N}S^{\frac{N}{2}}(\mu)<C_{MP}<\frac{2}{N}S^{\frac{N}{2}}(\mu)\leq\frac{1}{N}S^{\frac{N}{2}}\;  \hbox{if}\;\mu_2=\mu_1=\mu.$$
Combining     (\ref{bb-e5}), it follows that  both (\ref{hhe6-5}) and (\ref{hhe6-6}) are satisfied.
\ep
\vskip0.336in
\subsection{The case of $N=4, \frac{1}{2}<\big(\frac{1-\mu_1}{1-\mu_2}\big)^{\frac{3}{2}}, S^{\frac{N}{2}}(\mu_2)+S^{\frac{N}{2}}(\mu_1)\leq S^{\frac{N}{2}}$.}  In this case, we only consider $\alpha=\beta=2$.
If $\mu_2<\mu_1$, define
\be\lab{xue-e1}
M'_1:=1-\frac{4(1-\mu_1)^{\frac{3}{2}}}{[(1-\mu_2)^{\frac{3}{4}}+(1-\mu_1)^{\frac{3}{4}}]^2},
\ee
\be\lab{xue-e2}
M'_2:=2\big(\frac{1-\mu_1}{1-\mu_2}\big)^{\frac{3}{2}}-1,
\ee
\be\lab{xue-e3}
M'_3:=\frac{1}{2}\min\{M'_1,M'_2\},
\ee
then we have
\be\lab{xue-e4}
\frac{2}{N}(1-M'_3)\big(\frac{S(\mu_2)+S(\mu_1)}{2}\big)^{\frac{N}{2}}
>\frac{2}{N}S^{\frac{N}{2}}(\mu_1)>\frac{1+M'_3}{N}S^{\frac{N}{2}}(\mu_2).
\ee
If $\mu_2=\mu_1=\mu$, we define
\be\lab{xue-e5}
M'_3:=\frac{1}{4},
\ee
then we have
\be\lab{xue-e6}
\frac{2}{N}(1-M'_3)\big(\frac{S(\mu_1)+S(\mu_2)}{2}\big)^{\frac{N}{2}}>\frac{1+M'_3}{N}S^{\frac{N}{2}}(\mu_2).
\ee
Define
\be\lab{xue-e7}
M'_4:=\frac{1-\frac{1}{2}(1-M'_3)^{\frac{2}{N}}}{(\alpha+\beta)(\frac{1}{2})^{\frac{(N-2)(\alpha+\beta-2)}{4}}}.
\ee
If  $\alpha=\beta=2$, then  $M'_4=\frac{1}{2}-\frac{1}{4}(1-M'_3)^{\frac{1}{2}}>\frac{1}{4}$.
Similar to  Lemma \ref{lp-l1}, we have the following
\bl\lab{lp-l2}
Assume $N=4, \alpha=\beta= 2,  \frac{1}{2}<\big(\frac{1-\mu_1}{1-\mu_2}\big)^{\frac{3}{2}}$ such that  $\displaystyle S^{\frac{N}{2}}(\mu_1)+S^{\frac{N}{2}}(\mu_2)\leq S^{\frac{N}{2}}$.  Suppose that  the  weight function $h(x)$ satisfies $(H_2)$-$(H'_1)$ and that $\Theta<  C_{\alpha,\beta,\mu_1,\mu_2}$ (see  (\ref{zwm=may-2} ) ). If moreover,   $\Theta\leq  \min\{M'_4, C'_{\alpha,\beta}\}$,
 then $\overline{\Phi}$ has  a mountain pass geometry and the mountain pass level satisfies (\ref{amy-e-20}) and (\ref{amy-e-21}).
\el
\bp
Since we always assume $\alpha\geq 2, \beta\geq 2, \alpha+\beta\leq 2^*$, when $N=4$, the only possibility is that $\alpha=\beta=2, \alpha+\beta=4=2^*$. Thus we required $(H'_1)$. Note that  $\frac{1}{2}<\big(\frac{1-\mu_1}{1-\mu_2}\big)^{\frac{3}{2}}$ implies $S^{\frac{N}{2}}(\mu_2)<2S^{\frac{N}{2}}(\mu_1)$. We can follow carefully  the processes of Lemma \ref{lp-l1}   and obtain the  results. Here we omit   the details.
\ep
\vskip0.336in
\subsection{The case of $N=3, \frac{1}{2}<\frac{1-4\mu_1}{1-4\mu_2}, 2\big(\frac{S(\mu_2)+S(\mu_1)}{2}\big)^{\frac{N}{2}}> S^{\frac{N}{2}}$.}
In this case, if $\mu_2<\mu_1$, $M_3$ is redefined as
\be\lab{nn-e4}
M_3:=\frac{1}{2}\Big\{M_1, M_2, 1-\frac{\sqrt{2}}{[(1-4\mu_1)^{\frac{2}{3}}+(1-4\mu_2)^{\frac{2}{3}}]^{\frac{3}{2}}}\Big\},
\ee
where $M_1, M_2$ are defined in (\ref{xiu-e1}) and (\ref{xiu-e2}). Then (\ref{xiu-e4}) satisfies and
\be\lab{nn-e5}
2(1-M_3)\big(\frac{S(\mu_1)+S(\mu_2)}{2}\big)^{\frac{N}{2}}>S^{\frac{N}{2}}.
\ee
If $\mu_2=\mu_1=\mu$, $M_3$ is redefined as
\be\lab{nn-e6}
M_3:=\frac{1}{2}-\frac{1}{4(1-4\mu)},
\ee
then (\ref{lp-ez2}) and (\ref{nn-e5}) are satisfied. Here we    define
\be\lab{nn-e7}
M_4:=\frac{1-\frac{1}{2}(1-M_3)^{\frac{2}{N}}}{(\alpha+\beta)(\frac{1}{2})^{\frac{(N-2)(\alpha+\beta-2)}{4}}}>\frac{\sqrt{2}}{12}.
\ee
Then we have the following result:
\bl\lab{nn-lemma1}
Assume $N=3, \alpha\geq 2, \beta\geq 2, \alpha+\beta\leq 2^*,  \frac{1}{2}<\frac{1-4\mu_1}{1-4\mu_2}$  such that  $\displaystyle 2\big(\frac{S(\mu_2)+S(\mu_1)}{2}\big)^{\frac{N}{2}}> S^{\frac{N}{2}}$.   The  weight function $h(x)$ satisfies $(H_2)$ and $\begin{cases} (H_1)\;\hbox{if}\;\alpha+\beta<2^*\\ (H'_1)\;\hbox{if}\;\alpha+\beta=2^* \end{cases}$. In particular, when $\min\{\alpha, \beta\}=2$, we assume   that $\Theta<  C_{\alpha,\beta,\mu_1,\mu_2}$ (see  (\ref{zwm=may-2} ) ).    Then, if  $\Theta\leq  \min\{M_4, C'_{\alpha,\beta}\}$,  $\overline{\Phi}$ has  a mountain pass geometry and the mountain pass level satisfies  both (\ref{hhe6-5}) and (\ref{hhe6-6}).
\el
\bp Analogous to Lemma \ref{lp-l1}, $\overline{\Phi}$ exhibits a mountain pass geometry at  level $C_{MP}$ satisfying
$$\frac{2}{N}(1-M_3)\big(\frac{S(\mu_1)+S(\mu_2)}{2}\big)^{\frac{N}{2}}\leq C_{MP}<\frac{1}{N}S^{\frac{N}{2}}(\mu_1)+\frac{1}{N}S^{\frac{N}{2}}(\mu_2),$$
where $M_3$ is defined in (\ref{nn-e4}) for $\mu_2<\mu_1$ or (\ref{nn-e6}) for $\mu_2=\mu_1=\mu$.
Moreover, by (\ref{nn-e5}), we have
$$\frac{1}{N}S^{\frac{N}{2}}<C_{MP}<\frac{1}{N}S^{\frac{N}{2}}(\mu_1)+\frac{1}{N}S^{\frac{N}{2}}(\mu_2)\leq \frac{2}{N}S^{\frac{N}{2}}.$$
It follows that
$\frac{1}{N}S^{\frac{N}{2}}(\mu_2)<\frac{2}{N}S^{\frac{N}{2}}(\mu_1).$
If $\mu_2<\mu_1$, then we have $$\frac{1+M_3}{2}S^{\frac{N}{2}}(\mu_2)<\frac{2}{N}S^{\frac{N}{2}}(\mu_1)<C_{MP}<\frac{3}{N}S^{\frac{N}{2}}(\mu_1);$$
if $\mu_2=\mu_1=\mu$, we get that
$$\frac{1}{N}S^{\frac{N}{2}}(\mu_1)<\frac{1+M_3}{2}S^{\frac{N}{2}}(\mu_2)< C_{MP}<\frac{2}{N}S^{\frac{N}{2}}(\mu_1).$$
Thus,  $C_{MP}$ satisfies (\ref{hhe6-5}) and (\ref{hhe6-6}).
\ep
\vskip0.336in
\subsection{The case of $N=4, \frac{1}{2}<\big(\frac{1-\mu_1}{1-\mu_2}\big)^{\frac{3}{2}}, 2\big(\frac{S(\mu_2)+S(\mu_1)}{2}\big)^{\frac{N}{2}}> S^{\frac{N}{2}}$.}
In this case, if $\mu_2<\mu_1$, $M'_3$ is redefined as
\be\lab{ld-e1}
M'_3:=\frac{1}{2}\Big\{M'_1, M'_2, 1-\frac{2}{[(1-\mu_1)^{\frac{3}{4}}+(1-\mu_2)^{\frac{3}{4}}]^{2}}\Big\},
\ee
where $M'_1, M'_2$ are defined in (\ref{xue-e1}) and (\ref{xue-e2}). Then (\ref{xue-e4}) is satisfied  and
\be\lab{ld-e2}
2(1-M'_3)\big(\frac{S(\mu_1)+S(\mu_2)}{2}\big)^{\frac{N}{2}}>S^{\frac{N}{2}}.
\ee
If $\mu_2=\mu_1=\mu$, $M'_3$ is redefined as
\be\lab{ld-e3}
M'_3:=\frac{1}{2}-\frac{1}{4(1-4\mu)},
\ee
then (\ref{xue-e6}) and (\ref{ld-e2}) are satisfied.  Now $M'_4$ is   defined by
\be\lab{ld-e4}
M'_4:=\frac{1-\frac{1}{2}(1-M'_3)^{\frac{2}{N}}}{(\alpha+\beta)(\frac{1}{2})^{\frac{(N-2)(\alpha+\beta-2)}{4}}}.
\ee
 Similarly, we have the following result:
\bl\lab{nn-lemma2}
Assume $\displaystyle N=4, \alpha=\beta= 2,  \frac{1}{2}<\big(\frac{1-\mu_1}{1-\mu_2}\big)^{\frac{3}{2}}$ such that  $$\displaystyle 2\Big(\frac{S(\mu_2)+S(\mu_1)}{2}\Big)^{\frac{N}{2}}> S^{\frac{N}{2}}. $$   Suppose that the weight function $h(x)$ satisfies $(H_2)$ and $(H'_1)$;  $\Theta<  C_{\alpha,\beta,\mu_1,\mu_2}$ (see  (\ref{zwm=may-2})). Then, if $\Theta\leq  \min\{M'_4, C'_{\alpha,\beta}\}$,  $\overline{\Phi}$  has a mountain pass geometry at the  level $c$ satisfying  both (\ref{amy-e-20}) and (\ref{amy-e-21}).
\el
\bp  It is analogous to  the proof of  Lemmas \ref{lp-l2} and   \ref{nn-lemma1}.
\ep

Based  on the results of Lemma \ref{hl-6-1} $\sim$ Lemma \ref{nn-lemma2}, we can obtain the existence of mountain pass solution  to  problem (\ref{hhe6-4}).

\vskip0.336in
\noindent{\bf Proof of Theorem \ref{th-ling1}.}
Define
\begin{align*}
d_1:=&\min\Big\{0.3S^{\frac{1}{2}}(\mu_1),
10^{-3}\big[(1-4\mu_1)^{\frac{2}{3}}-\big(\frac{1}{2}\big)^{\frac{2}{3}}\big(1-4\mu_2\big)^{\frac{2}{3}}\big],5\times 10^{-4}\Big\}\\
=&\min\Big\{\frac{3}{10}, \frac{3}{10}S^{\frac{1}{2}}(\mu_1),\frac{\sqrt{2}}{12}, \frac{9}{100}, \frac{9}{100}S^{-\frac{1}{2}}(\mu_1),\\
& 10^{-3}\big[(1-4\mu_1)^{\frac{2}{3}}-\big(\frac{1}{2}\big)^{\frac{2}{3}}\big(1-4\mu_2\big)^{\frac{2}{3}}\big],5\times 10^{-4}\Big\}.
\end{align*}
Then $d_1<\min\{C_{\alpha, \beta, \mu_1,\mu_2}, M_4, C'_{\alpha,\beta},C_1, C_2\}. $  In particular,  the assumption (\ref{d111})$\Rightarrow \Theta\leq d_1$.
\begin{itemize}
\item Case I:\;$\displaystyle S^{\frac{N}{2}}(\mu_2)+S^{\frac{N}{2}}(\mu_1)\leq S^{\frac{N}{2}}$.
\end{itemize}
If $\min\{\alpha, \beta=2\}$,  then $\Theta<d_1$ implies that
$ \Theta<C_{\alpha, \beta, \mu_1, \mu_2}.$
Furthermore,
$$\Theta<d_1\Rightarrow \Theta<\min\{M_4, C'_{\alpha, \beta}\},$$
where $M_4$ is defined in (\ref{xiu-e5}). Then by Lemma \ref{lp-l1},  $\overline{\Phi}$ exhibits a mountain pass geometry and the mountain pass level satisfies (\ref{hhe6-5}) and (\ref{hhe6-6}).
\begin{itemize}
\item Case II:\;$\displaystyle 2\Big(\frac{S(\mu_1)+S(\mu_2)}{2}\Big)^{\frac{N}{2}}>S^{\frac{N}{2}}$.
\end{itemize}
For this case, $M_4$ is defined in (\ref{nn-e7}). Analogously,  by Lemma \ref{nn-lemma1}, $\overline{\Phi}$ has  a mountain pass geometry with energy  level satisfying   both (\ref{hhe6-5}) and (\ref{hhe6-6}).

For either case I or case II,   by the   Mountain Pass Theorem, there exists a sequence $\displaystyle \{(u_n, v_n)\}_{n\in \NN}\subset \overline{\mathcal{N}}$ such that
$$\overline{\Phi}(u_n, v_n)\rightarrow C_{MP},\; (\overline{\Phi})'\big|_{\overline{\mathcal{N}}}(u_n, v_n)\rightarrow 0\;\hbox{and}$$
$$\overline{\Phi}(u_n, v_n)>\frac{1+M_3}{N}S^{\frac{N}{2}}(\mu_2).$$
Recall that
$ \Theta<d_1$ implies that $ \Theta\leq \min\{C_1, C_2\},$
where $C_1, C_2$ are defined in (\ref{ling-e1}) and (\ref{ling-e2}).  By Lemma \ref{hl-6-1}, $\{(u_n, v_n)\}_{n\in \NN}$ admits a subsequence which converges strongly to a critical point $(u_0, v_0)$ of $\overline{\Phi}\big|_{\overline{\mathcal{N}}}$,
which is also a critical point of $\overline{\Phi}$ in $\D$. We observe that  $u_0\geq 0, v_0\geq 0, u_0v_0\not\equiv 0$, hence $(u_0, v_0)$ is a critical point of $\Phi$ in $\D$. That is,  $(u_0, v_0)$ is a nonnegative mountain pass solution of problem (\ref{hhe6-4}).\hfill$\Box$

\vskip0.336in
\noindent{\bf Proof of Theorem \ref{th-ling2}}.
Define
\begin{align*}
d_2:=&\min\Big\{\frac{\mu_1+\mu_2-(\mu_1+\mu_2)^2}{6}, 10^{-3}\big[(1-4\mu_1)^{\frac{2}{3}}-\big(\frac{1}{2}\big)^{\frac{2}{3}}\big(1-4\mu_2\big)^{\frac{2}{3}}\big],\\
&5\times 10^{-4}, \frac{3}{10}S^{\frac{1}{2}}(\mu_1)\Big\}\\
=&\min\Big\{\frac{\mu_1+\mu_2-(\mu_1+\mu_2)^2}{6}, 10^{-3}\big[(1-4\mu_1)^{\frac{2}{3}}-\big(\frac{1}{2}\big)^{\frac{2}{3}}\big(1-4\mu_2\big)^{\frac{2}{3}}\big],\\
&5\times 10^{-4}, \frac{3}{10}, \frac{3}{10}S^{\frac{1}{2}}(\mu_1),\frac{\sqrt{2}}{12},\frac{9}{100}, \frac{9}{100}S^{-\frac{1}{2}}(\mu_1)\Big\}.
\end{align*}
Then$$d_2< \min\Big\{\frac{1-(1-\varepsilon_1)^{\frac{2}{N}}}{2^{\frac{\beta}{2}}\alpha (1-\varepsilon_1)^{\frac{\alpha-2}{2^*}}}, \frac{1-(1-\varepsilon_1)^{\frac{2}{N}}}{2^{\frac{\alpha}{2}}\beta (1-\varepsilon_1)^{\frac{\beta-2}{2^*}}}, C_1, C_2,C_{\alpha, \beta, \mu_1,\mu_2}, M_4, C'_{\alpha,\beta}\Big\},
$$
and (\ref{d222})$\Rightarrow \Theta\leq d_2$.
Similar to the proof of Theorem \ref{th-ling1}, based  on the results of Lemma \ref{lp-l1} and Lemma \ref{nn-lemma1}, when $\Theta\leq d_2$, we obtain that
$\overline{\Phi}$ has a mountain pass geometry  which energy  level satisfies both (\ref{hhe6-5}) and (\ref{hhe6-6}).
By the  Mountain Pass Theorem, there exists a sequence $\displaystyle \{(u_n, v_n)\}_{n\in \NN}\subset \overline{\mathcal{N}}$ such that
$$\overline{\Phi}(u_n, v_n)\rightarrow C_{MP},\; (\overline{\Phi})'\big|_{\overline{\mathcal{N}}}(u_n, v_n)\rightarrow 0\;\hbox{and}$$
$$\overline{\Phi}(u_n, v_n)>\frac{1+M_3}{N}S^{\frac{N}{2}}(\mu_2).$$
Notice that
$$
\Theta\leq d_2\Rightarrow\Theta\leq \min\Big\{\frac{1-(1-\varepsilon_1)^{\frac{2}{N}}}{2^{\frac{\beta}{2}}\alpha (1-\varepsilon_1)^{\frac{\alpha-2}{2^*}}}, \frac{1-(1-\varepsilon_1)^{\frac{2}{N}}}{2^{\frac{\alpha}{2}}\beta (1-\varepsilon_1)^{\frac{\beta-2}{2^*}}}, C_1, C_2\Big\},$$
 where $\varepsilon_1$ satisfies (\ref{ml-e12}) and $C_1, C_2$ are defined in (\ref{ling-e1}) and (\ref{ling-e2}).
 By Lemma \ref{lemma-ml1}, $\{(u_n, v_n)\}_{n\in \NN}$ admits a subsequence which converges strongly to a critical point $(u_0, v_0)$ of $\overline{\Phi}\big|_{\overline{\mathcal{N}}}$, which is also a critical point of $\overline{\Phi}$ in $\D$. Also we know that $u_0\geq 0, v_0\geq 0, u_0v_0\not\equiv 0$, hence $(u_0, v_0)$ is a critical point of $\Phi$ in $\D$. That is,  $(u_0, v_0)$ is a nonnegative mountain pass solution of the problem (\ref{hhe6-4}).\hfill$\Box$

\vskip0.336in
\noindent{\bf Proof of Theorem \ref{th-ling3}.}
 Define
 \begin{align*}
 d_3:=&\min\Big\{\frac{2-(1-\mu_1)^{\frac{3}{2}}-(1-\mu_2)^{\frac{3}{2}}}{8}, \frac{2-\sqrt{2}}{4}\Big\}\\
 =&\min\Big\{\frac{2-(1-\mu_1)^{\frac{3}{2}}-(1-\mu_2)^{\frac{3}{2}}}{8}, \frac{2-\sqrt{2}}{4},\frac{1}{4},\frac{2-\sqrt{2}}{2}, \frac{\sqrt{2}}{4}\Big\}.
 \end{align*}
 Then $\displaystyle d_3 < \min\Big\{\frac{2-(1-\mu_1)^{\frac{3}{2}}-(1-\mu_2)^{\frac{3}{2}}}{8}, \frac{2-\sqrt{2}}{4}, M'_4, C'_{\alpha, \beta},C_{\alpha,\beta,\mu_1,\mu_2}\Big\}$ and (\ref{d333})$\Rightarrow \Theta\leq d_3$.  Thus by using  Lemma \ref{lemma-amy1}, Lemma \ref{lp-l2} and Lemma \ref{nn-lemma2},
the  problem (\ref{hsn-e1}) has a nontrivial weak solution
 $(u_0, v_0)$ such that $u_0\geq 0, v_0\geq0, u_0v_0\not\equiv 0$.\hfill$\Box$

\vskip0.23in






 \end{document}